\documentclass[10pt]{amsart}

\usepackage{amsmath}
\usepackage{amsfonts}
\usepackage{amssymb}
\usepackage{latexsym}
\usepackage{graphicx}
\usepackage{color}
\usepackage{fullpage}

\newtheorem{theorem}{Theorem}[section]
\numberwithin{equation}{section}

\def\p{u}
\def\q{y}
\def\dx{d\lambda}

\newcommand{\id} {{\rm Id}}

\renewcommand{\skew}{\mathop{\rm skew}}

\DeclareMathOperator{\sym}{sym}

\DeclareMathOperator{\Div}{Div}
\DeclareMathOperator{\Tr}{tr}
\DeclareMathOperator{\Curl}{Curl}

\DeclareMathOperator{\Grad}{Grad}

\DeclareMathOperator{\dev}{dev}
\DeclareMathOperator{\sL}{\mathfrak{sl}}
\DeclareMathOperator{\so}{\mathfrak{so}}

\newcommand{\Sym}{ {\rm{Sym}} }

\newcommand{\tr}[1]{ {\Tr \left({#1}\right)} }

\newcommand{\R}{\mathbb{R}}

\newcommand{\C}{\mathbb{C}}
\newcommand{\N}{\mathbb{N}}

\newcommand{\map}[3]
     {\begin{array}{ccccc}
     #1 & : & #2 & \longrightarrow & #3
     \end{array}}

\newcommand{\abb}[5]{#1 \hspace{1 mm} : \hspace{2 mm} #2
     \hspace{1 mm} \longrightarrow \hspace{1 mm} #3 \hspace{1 mm} ,
     \hspace{2 mm} #4 \hspace{1 mm} \longmapsto \hspace{1 mm} #5}

\newtheorem{definition}[theorem]{Definition}

\newtheorem{lemma}[theorem]{Lemma}

\begin{document}

\title{Dev-Div- and DevSym-DevCurl-Inequalities for Incompatible Square Tensor Fields with Mixed Boundary Conditions}
\author{Sebastian Bauer, Patrizio Neff, Dirk Pauly and Gerhard Starke}

\date{\today}

\maketitle

\vspace*{-5mm}
\begin{center}{\sf
Fakult\"at f\"ur Mathematik, Universit\"at Duisburg-Essen, Campus Essen\\
Thea-Leymann-Str.~9, D-45127 Essen, Germany
}\end{center}

\begin{abstract}\noindent
Let $\Omega\subset\R^n$, $n\geq2$, be a bounded Lipschitz domain and $1<q<\infty$.
We prove the inequality
$$\|T\|_{L^q(\Omega)}\leq C_{D\!D}\left(\|\dev T\|_{L^q(\Omega)}+\|\Div T\|_{L^q(\Omega)}\right)$$
being valid for tensor fields $T:\Omega\to\R^{n\times n}$  
with a normal boundary condition on some open and non-empty part $\Gamma_\nu$ of the
boundary $\partial\Omega$. Here $\dev T=T-\frac{1}{n}\tr{T}\cdot\id$ 
denotes the deviatoric part of the tensor $T$
and $\Div$ is the divergence row-wise. Furthermore, we prove
\begin{align*}
\|T\|_{L^2(\Omega)}
&\leq C_{D\!S\!C\!}\left(\|\dev\sym T\|_{L^2(\Omega)}+\|\Curl T\|_{L^2(\Omega)}\right)&
&\text{if }n\geq3,\\
\|T|\|_{L^2(\Omega)}
&\leq C_{D\!S\!D\!C\!}\left(\|\dev\sym T\|_{L^2(\Omega)}+\|\dev\Curl T\|_{L^2(\Omega)}\right)&
&\text{if }n=3,
\end{align*}
being valid for tensor fields $T$ with a tangential boundary condition on some 
open and non-empty part $\Gamma_\tau$ of $\partial\Omega$.
Here, $\sym T = \frac{1}{2}(T+T^\top)$ denotes the symmetric part of $T$ 
and $\Curl$ is the rotation row-wise.\\
{\bf Keywords:} 
Korn's inequality, Lie-algebra decomposition, 
Poincar\'e's inequality, Maxwell estimates, relaxed micromorphic model
\end{abstract}

\tableofcontents


\setcounter{equation}{0}

\section {Introduction}
`Every mathematical theorem has an inequality behind it ...' In this work
we consider $(n\times n)$-tensor fields $T$ on bounded domains $\Omega\subset \R^n$, $n\geq 2$,
with Lipschitz-continuous boundary $\partial \Omega$.
Such a tensor field may be decomposed pointwise orthogonally in its
{\sl symmetric part}  and its {\sl skew-symmetric part}
	\begin{equation}\label{symskew}
		T = {\rm sym\,}T+{\rm skew\,}T\,,
	\end{equation}
 where ${\rm sym\,}T = \frac{1}{2}(T+T^\top )$ and ${\rm skew\,}T=\frac{1}{2}\left(T-T^\top \right)$. 
In two recent papers \cite{neffpaulywitsch:13a,neffpaulywitsch:13b}, it has been shown that in $L^2(\Omega)$ 
the skew symmetric part of $T$  is controlled by the 
symmetric part and the ${\rm Curl}$ of $T$, leading to
	\begin{equation}\label{Sym-Curl-E}
		|\!| T |\!|_{L^2(\Omega)} \leq C_{S\!C}\left(|\!| {\rm sym\,}T |\!|_{L^2(\Omega)} +|\!| {\rm Curl\,}T |\!|_{L^2(\Omega)}	\right)\,,
	\end{equation}
if a {\sl tangential} boundary condition is imposed on some non-empty and open part $\Gamma_\tau$ of the boundary $\partial\Omega$. 
In classical terms $T\tau|_{\Gamma_\tau}=0$ is needed for all tangential-vectors $\tau$ on $\Gamma_\tau$. 
Here and hereafter all differential operators on tensor fields are taken row-wise.
For exact definitions of operators and function spaces, see Section 2.  
We shall call this inequality  the {\sl\bf Sym-Curl-inequality}.
 Since the  ${\rm Curl}$ operator vanishes on gradients, a certain variant of Korn's first inequality 
 follows immediately, i.e., 
with $T={\rm Grad\,}v$ and ${\rm Curl\,}{\rm Grad\,}=0$ we have 
	\begin{equation}\label{Korn-ineq-intro}
		|\!| {\rm Grad\,}v |\!|_{L^2(\Omega)} \leq C_{S\!C} \,|\!| {\rm sym\,}{\rm Grad\,}v |\!|_{L^2(\Omega)}
	\end{equation}
for all $v\in H^1(\Omega)$ with $ ({\rm Grad\,}v)\tau{\,|_{\Gamma_\tau}}=0$.
Obviously, this boundary condition is a weakening of the usual Dirichlet boundary condition $v\,|_{\Gamma_\tau}=0$, see the
discussion in \cite{neffpaulywitsch:12}.

The tensor $T$ may also be decomposed pointwise orthogonally in its 
{\sl trace-free} or {\sl deviatoric} and its  {\sl trace} or {\sl spherical part }
\begin{equation}\label{DevTr}
	T = {\rm dev\,}T+\frac{1}{n}{\rm tr}(T)\cdot{\rm Id}\,,
\end{equation}
where ${\rm Id}$ denotes the identity matrix in $\R^n$ and ${\rm tr}(T)=\sum_{i=1}^n T_{ii}$.

In Theorem \ref{Dev-Div-T} of this contribution, we show that in $L^q(\Omega)$, $1<q<\infty$, the trace part of $T$ 
is controlled by the deviatoric part and the divergence of $T$, i.e.,
\begin{equation}\nonumber
	|\!| T |\!|_{L^q(\Omega)} \leq C_{D\!D}\left(|\!| {\rm dev\, }T |\!|_{L^q(\Omega)}+|\!| {\rm Div\,} T|\!|_{L^q(\Omega)}\right)\,,
\end{equation}
if a {\sl normal} boundary condition is imposed  on some non-empty and open  part  $\Gamma_\nu$ of $\partial\Omega$. 
In classical terms
\begin{equation}\label{boundaryconditionDevDivclassicalversion}
		T \nu{\,|_{\Gamma_\nu}}=0
	\end{equation} 
is needed for the normal vector $\nu$ at $\Gamma_\nu$.
We shall call this inequality the {\bf Dev-Div-inequality}.

In case that $n=3$ and $T={\rm Curl\,}S$ is already a ${\rm Curl\,}$ 
of a tensor field $S$ having the proper tangential boundary condition on $\Gamma_\tau$,
we conclude that $T$ is already controlled by its deviatoric part alone, i.e.,
\begin{equation}\label{DevDiv-estimate}
	\begin{array}{rcl}
		|\!| {\rm Curl\,}S |\!|_{L^q(\Omega)} 	& \leq &  C_{D\!D}|\!| {\rm dev\, }{\rm Curl\,}S|\!|_{L^q(\Omega)}\,,
	\end{array}
\end{equation}
since ${\rm Div\,}{\rm Curl\,}=0$ and $T$ inherits the proper normal boundary condition from $S$.
The inequality \eqref{DevDiv-estimate} may be seen as a Korn-type inequality, cf. \eqref{Korn-ineq-intro}.
Both orthogonal decompositions (\ref{symskew}) and (\ref{DevTr}) may be combined by appealing to the Cartan-decomposition
of the Lie-algebra $\mathfrak{gl}(n)$
	\begin{equation}\nonumber\label{DevSymSkewTr}
		\begin{array}{ccccccc}
			\mathfrak{gl}(n) & = & \left(\mathfrak{sl}(n)\cap\Sym(n)\right)  & \oplus &  \mathfrak{so}(n) & \oplus &  \R\cdot {\rm Id}\, \\[1.5ex]
			T 		 & = & {\rm dev\,}{\rm sym\,}T   &    +      & {\rm skew\,}T &+&\displaystyle\frac{1}{n}{\rm tr}(T)\cdot {\rm Id}\,.
		\end{array}
	\end{equation}
Here, $\mathfrak{sl}(n)$ denotes the Lie-algebra of trace free matrices and
$\mathfrak{so}(n)$ denotes the Lie-algebra of skew-symmetric matrices in $\R^{n\times n}$.
Now, in a naive manner an estimate of the following kind could be guessed 
	\begin{equation}\nonumber
		|\!| T |\!|_{L^2(\Omega)} \leq C \left(|\!| {\rm dev\,}{\rm sym\,}T |\!|_{L^2(\Omega)}
				+|\!| {\rm Curl\,}T |\!|_{L^2(\Omega)}+|\!| {\rm Div\,}T |\!|_{L^2(\Omega)}\right),
	\end{equation}
accompanied by suitable boundary conditions.
In fact, in Theorem \ref{DevSym-DevCurl-T} we prove a somewhat stronger result: For $n=3$ we prove the new 
{\bf DevSym-DevCurl-inequality}
	\begin{equation}\nonumber
		|\!| T |\!|_{L^2(\Omega)} \leq C_{D\!S\!D\!C}\left(|\!| {\rm dev\,}{\rm sym\,}T |\!|_{L^2(\Omega)}
		+|\!| {\rm dev\,}{\rm Curl\,}T |\!|_{L^2(\Omega)}\right)\,,
	\end{equation}
where again a tangential boundary condition is imposed on some
non-empty and open part $\Gamma_\tau$ of the boundary.
Since the deviatoric part is only defined for quadratic tensors, this estimate does not make sense for 
 $n\not=3$ . In general,
${\rm Curl\,} T$ is a $(n(n-1)/2\times n)$-matrix and we prove that
for $n\geq3$
	\begin{equation}\label{DevSymCurl-intro}
		|\!| T |\!|_{L^2(\Omega)} \leq C_{D\!S\!C}\left(|\!| {\rm dev\,}{\rm sym\,}T |\!|_{L^2(\Omega)}
		+|\!| {\rm Curl\,}T |\!|_{L^2(\Omega)}\right)
	\end{equation}
holds. In order to show \eqref{DevSymCurl-intro} we first prove for $n\geq 3$ a Korn type inequality, i.e.
	\begin{equation}\label{DevSymGrad-E}
		|\!| {\rm Grad\,}v |\!|_{L^2(\Omega)} \leq C\, |\!| {\rm dev\,}{\rm sym\,}{\rm Grad\,}v |\!|_{L^2(\Omega)}	
	\end{equation}
for all $v\in H^1(\Omega)$ with $({\rm Grad\,} v)\tau|_{\Gamma_\tau}=0$. 
The results of
this paper and some of their consequences have already been sketched in \cite{bauerneffpaulystarke:13a}
and \cite{bauerneffpaulystarke:13b}.

Whereas inequalities of  Sym-Curl-type are investigated by some of the present authors in a series of papers
for the first time, see \cite{neffpaulywitsch:11a, 
neffpaulywitsch:11b, neffpaulywitsch:12, neffpaulywitsch:13a, neffpaulywitsch:13b}, 
 there are already several contributions  to Div-Dev-type inequalities in the literature:  
In  \cite[Lemma 3.1]{arnolddouglasgupta:84} 
a Div-Dev-estimate is proved for  $n=2$ replacing the boundary condition by the average condition 
$\int_\Omega {\rm tr}(T){\,\rm d}x=0$. 
In the proof  of Theorem  \ref{Dev-Div-T} we adopt the idea of proof from this Lemma. 

In \cite[Lemma 3.2]{caistarke:03} for $n=2$ and $n=3$ the estimate
	\begin{equation}\nonumber
		|\!| T |\!|_{L^2(\Omega)}^2 \leq C \left(\frac{1}{2\mu}|\!| {\rm dev\,}T |\!|_{L^2(\Omega)}^2+
			\frac{1}{n(n\lambda+2\mu)}|\!| {\rm tr}(T) |\!|_{L^2(\Omega)}^2+|\!| {\rm Div\,}T |\!|^2_{H^{-1}(\Gamma_\tau;\,\Omega)}\right)
\end{equation}
is shown by means of a Helmholtz decomposition.
In  the notation used in this paper  $H^{-1}(\Gamma_\tau;\,\Omega)$ denotes 
the dual space of $H\!\left({\rm Grad};\,\Gamma_\tau;\,\Omega\right)$.
This estimate holds uniformly in $0<\mu_1\leq \mu\leq \mu_2$ and $0<\lambda<\infty$. Therefore, in the  (incompressible) 
limit $\lambda\to
\infty$ this estimate implies a Dev-Div-estimate. 

Korn-type estimates, replacing the symmetric gradient by its trace-free part are given in \cite{dain:06} and
\cite{Neff_JeongMMS08}, i.e.
\begin{equation}\label{KornSecondTraceFree-E}
	|\!| {\rm Grad\,} v |\!|_{L^2(\Omega)}\leq C\left(|\!| {\rm dev\,}{\rm sym\;}{\rm Grad\,}v |\!|_{L^2(\Omega)}
	+|\!| v |\!|_{L^2(\Omega)}\right)
\end{equation}
for all $v \in H\!\left({\rm Grad};\,\Omega\right)$ and $n\geq 3$. 
In \cite[Theorem 3.2]{Reshetnyak94} a trace-free version of Korn's first inequality is shown by means of integral representations.
In detail it is shown that for $1<q<\infty$ and any projector $\Pi$ from $W^q\left({\rm Grad};\,\Omega\right)$ onto the finite dimensional
kernel of ${\rm dev\,}{\rm sym\,}{\rm Grad}$, there exists a constant $C>0$, such that for all $u\in W^q\left({\rm Grad};\,\Omega\right)$
  \begin{equation}\nonumber
		|| u -\Pi u||_{W^q\left({\rm Grad};\,\Omega\right)}\leq C\,|| {\rm dev\,}{\rm sym\,}{\rm Grad\,}u ||_{L^q(\Omega)}\,.
  \end{equation}

It is well known, that for $n=2$ estimate (\ref{KornSecondTraceFree-E})
fails to hold true, since in this case the kernel of ${\rm dev\,}{\rm sym\;}{\rm Grad\,}$ is given by the holomorphic 
functions and thus is infinite-dimensional. On the other hand, in \cite[Appendix]{neffpaulywitsch:13a} inequality (\ref{DevSymGrad-E}) 
is proved for $v\in H\!\left({\rm Grad};\,\partial\Omega;\,\Omega\right)$ by simple partial integration and 
some elementary estimates. In \cite{fuchsschirra:09} it is proved that 
$$ |\!| {\rm Grad\,}v |\!|_{L^q(\Omega)} \leq C |\!| {\rm dev\,}{\rm sym\,}{\rm Grad\,}v |\!|_{L^q(\Omega)}$$ holds
for $v\in W^q\!\left({\rm Grad};\,\partial\Omega;\,\Omega\right)$ for $n=2$ and $1<q<\infty$, and in \cite{fuchsrepin:10}
this inequality is proved for $q=1$, $v\in W^1\!\left({\rm Grad};\,\partial\Omega;\,\Omega\right)$ 
and arbitrary space dimensions $n$.
In  Section 6 we show that for the case of only a partial boundary condition, i.e.
 $v\in H\!\left({\rm Grad};\,\Gamma_\tau;\,\Omega\right)$,
the estimate  (\ref{DevSymGrad-E}) is false by means of a construction taken from \cite{pompe:10}.

What about  inequalities like DevSym-DevSymCurl or other combinations? In Section 6 we give some negative 
results in that direction.  It may be quite illuminating to see by some
{\sl simple} arguments, why the kernel of the operators defining the right hand side of our 
inequalities are trivial on, say, the space of smooth compactly supported tensor fields.
 Some calculations in that direction are also presented in Section 6.
In Section 7  applications of the derived inequalities are given.
 The remaining part of the paper is organized as follows: In Section 2 we shall give notations and definitions 
used in this paper. In Section 3 we provide the proof of the Dev-Div-inequality
and in Section 4 and 5 we give the proofs of the DevSym-Curl- and the DevSym-DevCurl-inequality. 
In the Appendix we prove a representation formula
for the kernel of ${\rm dev\,}{\rm sym\;}{\rm Grad\,}$ in arbitrary space dimensions 
used in the proof of Theorem \ref{DevSym-DevCurl-T}.

\section{Definitions and preliminaries}

Throughout the entire paper we assume $\Omega\subset\R^n$, $n\geq2$, to be a bounded domain with  boundary $\partial\Omega$.
Moreover, let $\Gamma_\tau$ be a relatively open subset of $\partial\Omega$ and $\Gamma_\nu:=\partial\Omega\setminus \bar{\Gamma}_\tau$.
Here, the subscripts $\tau$ and $\nu$ refer to the {\sl tangential} and {\sl normal} boundary condition,
respectively. 

The usual Lebesgue-spaces of $q$-integrable functions, vector fields and tensor fields on $\Omega$
with values in $\R$, $\R^n$ and $\R^{n\times n}$, respectively, will be denoted by $L^q(\Omega)$.
Moreover, we introduce the standard Sobolev-spaces
	\begin{eqnarray}\label{DefinitionWgrad}\nonumber
			W^q\!\left({\rm grad};\,\Omega\right) & := & 
				\left\{u \in L^q(\Omega)\; | \; {\rm grad\;}u \in L^q(\Omega) \right\}=W^{1, q}(\Omega)\,,\\
				\label{DefinitionWdiv}\nonumber
			W^q\!\left({\rm div};\,\Omega\right)	& := &	\left\{u \in L^q(\Omega)\; | \; {\rm div\;}u \in L^q(\Omega) \right\}\,,\\
				\label{DefinitionWcurl}\nonumber
			W^q\!\left({\rm curl};\,\Omega\right)	& := &	\left\{u \in L^q(\Omega)\; | \; {\rm curl\;}u \in L^q(\Omega) \right\}\,,
	\end{eqnarray}
where ${\rm grad}$, ${\rm div}$ and ${\rm curl\,}$ are the usual differential operators gradient, divergence and 
rotation\footnote{For a definition of the rotation for $n\not=3$, see, e.g. \cite{neffpaulywitsch:13b}.}, respectively.
All derivatives are understood in the distributional sense. For $q=2$ we replace as usual $W^2$ by $H$. 

In order to realize certain boundary conditions we make use of the spaces
	\begin{equation}\label{Definitionsmoothfuntions}\nonumber
		C^\infty(\Gamma,\, \bar{\Omega}) := \left\{u{|_\Omega}\; | \; u \in C^\infty_0(\R^n\setminus \bar{\Gamma})\right\}
	\end{equation}
	for $\Gamma= \partial\Omega, \Gamma_\tau$ or $\Gamma_\nu$
	and define 
	\begin{equation}\label{Definitionboundaryconditions}
			W^q\!\left({\rm grad};\,\Gamma_\tau;\,\Omega\right), \quad 
			W^q\!\left({\rm div};\,\Gamma_\nu;\,\Omega\right)\quad\text{and}\quad
			W^q\!\left({\rm curl};\,\Gamma_\tau;\,\Omega\right)
	\end{equation}
as completion under the respective graph norms of the scalar-valued space $C^\infty(\Gamma_\tau,\, \bar{\Omega})$ and the 
vector-valued spaces $C^\infty(\Gamma_\nu,\, \bar{\Omega})$ and $C^\infty(\Gamma_\tau,\, \bar{\Omega})$, respectively.
Therefore, these spaces generalize the homogeneous Dirichlet boundary conditions
	\begin{equation}\nonumber
		u{|_{\Gamma_\tau}} = 0 \,\text{(scalar),}\qquad
		\nu\cdot v{|_{\Gamma_\nu}} = 0\,\text{(normal)}\qquad\text{and}\qquad
		\nu\times v{|_{\Gamma_\tau}} = 0\,\text{(tangential), }
	\end{equation}
respectively. 

Now we extend our notations to vector and tensor fields by defining all differential operations on rows.
Thus, for a vector field
$v=(v_1, \dots, v_n)^\top $ we define the tensor field 
${\rm Grad\,}v:=({\rm grad}^\top v_1,\dots, {\rm grad}^\top v_n)^\top $, where ${}^\top $ denotes the transpose. 
Note, that ${\rm Grad\,}v$ is  just the Jacobian of $v$.
For a tensor field $T$ we define the divergence ${\rm Div\,}T:=\left({\rm div\,}T_1^\top , \dots, {\rm div\,}T_n^\top \right)^\top $ and the rotation
$ {\rm Curl\,}T= \left({\rm curl}^\top T_1^\top ,\dots, {\rm curl}^\top T_n^\top \right)^\top $, where
$T_i$ denote the row-vectors of $T$, i.e., $T=\left(T_1,\dots,T_n\right)^\top $.
The corresponding Sobolev-spaces will be denoted by 
$$W^q\!\left({\rm Grad};\,\Omega\right),\quad 
H\!\left({\rm Grad};\,\Omega\right),\quad 
W^q\!\left({\rm Grad};\,\Gamma_\tau;\,\Omega\right),\quad 
H\!\left({\rm Grad};\,\Gamma_\tau;\,\Omega\right)$$ 
and so on.
Note that the spaces $W^q\!\left({\rm Div};\,\Gamma_\nu;\,\Omega\right)$ and  $H\!\left({\rm Div};\,\Gamma_\nu;\,\Omega\right)$
generalize the {\it normal boundary condition} $T\nu|_{\Gamma_\nu}=0$, while the spaces
$W^q\!\left({\rm Curl};\,\Gamma_\tau;\,\Omega\right)$ and $H\!\left({\rm Curl};\,\Gamma_\tau;\,\Omega\right)$
generalize the {\sl tangential boundary condition} $T\tau|_{\Gamma_\tau}=0$.

In general, we only assume fairly weak regularity assumptions on the boundary. To be specific, 
from the theory of scalar valued functions we need the compact embedding
of  $W^{1, q}(\Omega)$ into $L^q(\Omega)$, i.e. Rellich's selection theorem,
Korn's second inequality in $L^q(\Omega)$ and the so-called Lions-Lemma (\ref{LionsLemma-q}),
which are guaranteed, if the boundary $\partial\Omega$ is locally the graph of a Lipschitz-continuous function, see e.g. \cite{adams:75} and \cite{necas:67b}. Moreover, from the theory of vector fields, we need the so-called
Maxwell compactness property for mixed boundary conditions, i.e., the compact embedding of 
$H\!\left({\rm curl};\,\Gamma_\tau;\,\Omega\right)\cap H\!\left({\rm div};\,\Gamma_\nu;\,\Omega\right)$
into $L^2(\Omega)$. This implies also for tensor fields the Maxwell estimate \eqref{Maxwell-E}
and the Helmholtz decomposition \eqref{HelmholtzDecomposition}, 
which are also essential tools in our arguments. 
These hold for Lipschitz boundaries $\partial\Omega$ as well, provided that the interface
$\bar{\Gamma}_\tau\cap\bar{\Gamma}_\nu$ is Lipschitz itself.
Therefore, throughout this paper we will assume generally the latter regularity.

\section{The Dev-Div-inequality}

In this section we shall prove the following theorem.
\begin{theorem}\label{Dev-Div-T}
	Let $\Gamma_\nu\not=\emptyset$ and $1<q<\infty$. 
	Then there exists a constant $C_{D\!D}$, 
	such that 
	the following estimates hold:
	\begin{enumerate}
		\item[(i)] For all $T\in W^q\!\left({\rm Div};\,\Gamma_\nu;\,\Omega\right)$
		\begin{equation}\label{Dev-Div-E}
			|\!| T |\!|_{L^q(\Omega)} \leq C_{D\!D}\left(|\!| {\rm dev\,}T |\!|_{L^q(\Omega)}+	|\!| {\rm Div\,}T |\!|_{L^q(\Omega)}\right)
		\end{equation}
		\item[(i')] 		and 
		\begin{equation*}
			|\!| T |\!|_{W^q\!\left({\rm Div};\,\Omega\right)} \leq C_{D\!D}\left(|\!| {\rm dev\,}T |\!|_{L^q(\Omega)}+	|\!| {\rm Div\,}T |\!|_{L^q(\Omega)}\right).
		\end{equation*}
		\item[(ii)] If $n=3$, for all $T\in W^q\!\left({\rm Curl};\,\Gamma_\nu;\,\Omega\right)$
		\begin{equation}\label{Curl-DevCurl-E}
			|\!| {\rm Curl\,} T |\!|_{L^q(\Omega)} \leq C_{D\!D}\,|\!| {\rm dev\,}{\rm Curl\,} T |\!|_{L^q(\Omega)}\,.
		\end{equation}
	\end{enumerate}
\end{theorem}
We shall prove this Theorem using the following lemma, guaranteeing the existence of some suitable divergence-potential.
\begin{lemma}\label{DivergenceLemma}
	Let $\Gamma_\nu\not=\emptyset$ and $1<q<\infty$. Then, there exists a constant $C>0$, 
such that for all real-valued functions 
	$g\in L^q(\Omega)$ there is a vector field $v\in W^q\!\left({\rm Grad};\,\Gamma_\tau;\,\Omega\right)$ with
		\begin{equation}\label{LemmaDivergence}
			{\rm div\,}v = g \quad\text{and}\quad
			|\!| v |\!|_{W^q\!\left({\rm Grad};\,\Omega\right)} \leq C \,|\!| g |\!|_{L^q(\Omega)}\,.
		\end{equation}
\end{lemma}
In the case $\Gamma_\nu=\emptyset$, this Lemma has been proved in \cite[Lemma 2.1.1]{sohr:01}  under 
the additional normalization assumption $\int_\Omega g\,\dx=0$. 
With minor modifications the same proof  also works in the situation
under consideration. For the convenience of the reader we shall give it in some detail.\\[2ex]
{\sl Proof:} 
	The linear operator 
		\begin{equation}\nonumber
			\abb{\underline{{\rm div}}}{W^q\!\left({\rm Grad};\,\Gamma_\tau;\,\Omega\right)}{L^q(\Omega)}{v}{{\rm div\,}v}
		\end{equation}
	is bounded, i.e.
		\begin{equation}\label{div-e}
			|\!| \underline{\rm div\,}v |\!|_{L^q(\Omega)}\leq C_1|\!| v |\!|_{W^q\!\left({\rm Grad};\,\Omega\right)}
		\end{equation}
	holds for all $v\in W^q\!\left({\rm Grad};\,\Gamma_\tau;\,\Omega\right)$.
	Now, we identify $L^q(\Omega)^\prime= L^{q'}(\Omega)$, where 
$1/q+1/q'=1$. 
	Further, we consider the dual operator of $\underline{{\rm div}}$, 
		\begin{equation}\nonumber
			\map{\underline{{\rm div}}^\prime =-{\rm g\underline{rad}}\,}{L^{q'}(\Omega)}
			{W^q\!\left({\rm Grad};\,\Gamma_\tau;\,\Omega\right)^\prime=:W^{-1,q'}({\rm Grad};\,\Gamma_\nu;\,\Omega)}\,,
		\end{equation}
	defined by
		\begin{equation}\nonumber
			-\langle{\rm g\underline{rad}}\, \,\p, v\rangle \,:=\, \int_\Omega \p\, {\rm div\,}v\,\dx
		\end{equation}
	for all $v\in W^q\!\left({\rm Grad};\,\Gamma_\tau;\,\Omega\right)$ and all $\p\in L^{q'}(\Omega)$.
	Here $\langle\,\cdot\, , \,\cdot\, \rangle$ denotes the duality pairing of 
	$W^q\!\left({\rm Grad};\,\Gamma_\tau;\,\Omega\right)$ and $W^{-1,q'}({\rm Grad};\,\Gamma_\nu;\,\Omega)$.
	Utilizing (\ref{div-e}) and the definition of the norm in the dual space we obtain
	the continuity of ${\rm g\underline{rad}}\,$, i.e.,
		\begin{equation}\nonumber
			|\!| {\rm g\underline{rad}}\, \p |\!|_{W^{-1,q'}({\rm Grad};\,\Gamma_\nu;\,\Omega)} \leq C_1 \,|\!| \p |\!|_{L^{q'}(\Omega)}\,.
		\end{equation}
	We will show that also the reversed inequality holds true: 
	There exists a constant $C_2>0$, 
	such that for all $\p\in L^{q'}(\Omega)$
		\begin{equation}\label{Inequ1}
			|\!| \p |\!|_{L^{q'}(\Omega)} \leq C_2 \,|\!| {\rm g\underline{rad}}\, \p |\!|_{W^{-1,q'}({\rm Grad};\,\Gamma_\nu;\,\Omega)}\,.
		\end{equation}
	To prove  (\ref{Inequ1}) we use the usual contradiction argument: Assume the inequality  is false, then there exists a sequence
	$(\p_j)\subset L^{q'}(\Omega)$ with
		\begin{equation}\label{Annahme}
			 |\!| \p_j |\!|_{L^{q'}(\Omega)}=1
			\;\text{for all $j$}\qquad\text{and}\qquad
			\lim_{j\to \infty}|\!| {\rm g\underline{rad}}\,\p |\!|_{W^{-1,q'}({\rm Grad};\,\Gamma_\nu;\,\Omega)}=0\,.
		\end{equation}
	Since $(\p_j)$ is bounded in $L^{q'}(\Omega)$, by weak compactness there exists a subsequence of $(\p_j)$, 
	also called $(\p_j)$, 
	and a  $\p\in L^{q'}(\Omega)$, such that 
		\begin{equation}\nonumber
			\p_j \rightharpoonup \p \quad\text{weakly in}\quad L^{q'}(\Omega)\,.
		\end{equation}
	Since for all $v\in W^q\!\left({\rm Grad};\,\Gamma_\tau;\,\Omega\right)$
		\begin{eqnarray}\label{Null}
			|\langle{\rm g\underline{rad}}\, \p, v\rangle| & = & \left|\int_\Omega \p\,{\rm div\,}v\,\dx\right| 
										   =  	\lim_{j\to\infty}\left|\int_\Omega \p_j\,{\rm div\,}v\,\dx\right| \\ \nonumber
										&   =  & \lim_{j\to\infty}|\langle{\rm g\underline{rad}}\, \p_j, v\rangle| 
										 \leq 
										   \lim_{j\to\infty}|\!| {\rm g\underline{rad}}\, \p_j |\!|_{W^{-1,q'}({\rm Grad};\,\Gamma_\nu;\,\Omega)}\,
													|\!| v |\!|_{W^q\!({\rm Grad\,}\Omega)}=0,
		\end{eqnarray}
	  we conclude ${\rm g\underline{rad}}\, \p=0$, which implies
	 ${\rm grad\,} \p=0$ in the distributional sense and hence by the fundamental lemma $\p=const$,  
	 see also e.g. \cite[II, (1.7.18)]{sohr:01}.
	As $\Gamma_\nu\not=\emptyset$ is relatively open, 
	there exists a vector field $v \in W^q\!\left({\rm Grad};\,\Gamma_\tau;\,\Omega\right)$ such that 
		\begin{equation}\nonumber
			\int_{\Omega}{\rm div\,} v \,\dx \not=0\,.
		\end{equation}
	Employing this, (\ref{Null}) and $\p=const$  we conclude $\p = 0$.	
	Remarkably, the operator ${\rm g\underline{rad}}\,$, although being a kind of differential operator, does not vanish on constant 
	functions.
	
	Following \cite[Theorem 6.3]{adams:75} the embedding 
	$W^q\!\left({\rm grad};\,\Omega\right)\hookrightarrow  L^{q}(\Omega)$ is compact.
	Hence, of course also $W^q\!\left({\rm grad};\,\Gamma_\tau;\,\Omega\right)\hookrightarrow  L^{q}(\Omega)$ is compact.
	Using \cite[X.4]{yoshida:80}, the dual embedding 
	$$L^{q'}(\Omega)\hookrightarrow W^{-1,q'}({\rm grad};\,\Gamma_\nu;\,\Omega)
	:=W^q\!\left({\rm grad};\,\Gamma_\tau;\,\Omega\right)',$$
	defined by $\langle f, w\rangle=\int_\Omega f\cdot w \,\dx$ for all $f\in L^{q'}(\Omega)$ and 
	$w\in W^q\!\left({\rm grad};\,\Gamma_\tau;\,\Omega\right)$,
	is compact as well.
	Thus, we can select a subsequence, again denoted by $(\p_j)$, which converges to some 
	$\hat{\p}\in W^{-1,q'}({\rm grad};\,\Gamma_\nu;\,\Omega)$ in 
	$W^{-1,q'}({\rm grad};\,\Gamma_\nu;\,\Omega)$.
	As we have seen, $(\p_j)$ also converges weakly in  $L^{q'}(\Omega)$ to $\p=0$ and therefore we get  $\hat{\p}=0$.
	Now  we use the so-called  {\sl Lions-Lemma} from \cite{necas:67b} (concerning the history of the Lions-Lemma, see also 
	\cite{Ciarlet10}):
	There is a positive constant $C_3$, such that for all $u\in L^{q'}(\Omega)$
		\begin{equation}\label{LionsLemma-q}
			|\!| u |\!|_{L^{q'}(\Omega)} \leq C_{3}\left(|\!| {\rm grad\,} u |\!|_{W^{-1,q'}({\rm Grad};\,\Omega)}
			+|\!| u |\!|_{W^{-1,q'}({\rm grad};\,\Omega)}\right)\,,
		\end{equation}
	where we set
	$W^{-1,q'}({\rm grad};\,\Omega):=W^{q}\!\left({\rm grad};\,\partial\Omega;\,\Omega\right)'$
	and $W^{-1,q'}({\rm Grad};\,\Omega):=W^{q}\!\left({\rm Grad};\,\partial\Omega;\,\Omega\right)'$. 
	The norms of  dual spaces
	$W^{-1,q'}({\rm grad};\,\Gamma_\nu;\,\Omega)$ and $W^{-1,q'}({\rm Grad};\,\Gamma_\nu;\,\Omega)$
	 are stronger than the norms of 
	 $W^{-1,q'}({\rm grad};\,\Omega)$ and $W^{-1,q'}({\rm Grad};\,\Omega)$. Hence, we can estimate
		\begin{eqnarray}\nonumber
			1=|\!| \p_j |\!|_{L^{q'}(\Omega)} 	& \leq & C_3\left(|\!| {\rm grad\,} \p_j |\!|_{W^{-1,q'}({\rm Grad};\,\Omega)}
												+|\!| \p_j |\!|_{W^{-1,q'}({\rm grad};\,\Omega)}\right) \\ \nonumber
										& \leq & C_3\left(|\!| {\rm g\underline{rad}}\, \p_j |\!|_{W^{-1,q'}({\rm Grad};\,\Gamma_\nu;\,\Omega)}
														+|\!| \p_j |\!|_{W^{-1,q'}({\rm grad};\,\Gamma_\nu;\,\Omega)}\right)
													\longrightarrow 0
		\end{eqnarray}
	for  $j\to\infty$, in contradiction to  (\ref{Annahme}). Thus (\ref{Inequ1}) is proved. 
	
  By (\ref{Inequ1}), 
	the range $R({\rm g\underline{rad}}$) of the operator ${\rm g\underline{rad}}$ is a closed subspace of
	$W^{-1,q'}({\rm Grad};\,\Gamma_\nu;\,\Omega)$.  Since $R( {\rm g\underline{rad}})$ is the range of the dual operator of 
	 $\underline{{\rm div}}$, the {\sl closed range theorem}, see e.g. \cite[VII.5]{yoshida:80},  
	yields that the range $R(\underline{{\rm div}})$ is  also closed and we have
		\begin{equation}\nonumber
			R(\underline{{\rm div}})=\left\{f\in L^q(\Omega) \; : \; \int_\Omega f\cdot \p\,\dx =0
									\quad\text{for all} \quad \p\in N({\rm g\underline{rad}})\right\}\;,
		\end{equation}
	where $N({\rm g\underline{rad}})$ denotes the kernel of the operator ${\rm g\underline{rad}}$.
	We have already shown above that ${\rm g\underline{rad}}\,\p=0$ implies $\p= 0$, i.e.
	$N({\rm g\underline{rad}})=\{0\}$.
	Therefore, 
		\begin{equation}\label{rangediv}
			R(\underline{{\rm div}}) = L^q(\Omega)\,.
		\end{equation}
	In order to get the estimate in (\ref{LemmaDivergence}), we consider the quotient space
		\begin{equation}\nonumber
			W^q\!\left({\rm Grad};\,\Gamma_\tau;\,\Omega\right)/N(\underline{{\rm div}}) :=
				\left\{[v] \; | \; v\in W^q\!\left({\rm Grad};\,\Gamma_\tau;\,\Omega\right)\right\}\,,
		\end{equation}
	with $[v]:=v+N(\underline{{\rm div}})$, $v\in W^q\!\left({\rm Grad};\,\Gamma_\tau;\,\Omega\right)$  and the associated norm
		\begin{equation}\nonumber
			|\!| [v] |\!|_{W^q\!\left({\rm Grad};\,\Gamma_\tau;\,\Omega\right)/N(\underline{\rm div})} 
			:=\inf_{w\in N(\underline{{\rm div}})}|\!|v+w|\!|_{W^q\!\left({\rm Grad};\,\Omega\right)}\,.
		\end{equation}
	Thus, the linear operator
		\begin{equation}\nonumber
			\abb{\overline{{\rm div}}}{W^q\!\left({\rm Grad};\,\Gamma_\tau;\,\Omega\right)/N(\underline{{\rm div}})}
			{L^q(\Omega)}{\qquad[v]}{{\rm div\,}v}
		\end{equation}
	is well-defined, bijective and bounded. According to the bounded inverse theorem, see e.g. \cite[II.5]{yoshida:80},
	the inverse operator $\overline{{\rm div}}^{\,-1}$, mapping $L^q(\Omega)$ to  
	$W^q\!\left({\rm Grad};\,\Gamma_\tau;\,\Omega\right)/N(\underline{{\rm div}})$ is bounded.
	Hence there exists a constant $C_4>0$, such that for all
	$g\in L^q(\Omega)$ with $g={\rm div\,}v$ and $v\in W^q\!\left({\rm Grad};\,\Gamma_\tau;\,\Omega\right)$
	 	\begin{equation}\nonumber
				\inf_{w\in N(\underline{{\rm div}})}|\!|v+w|\!|_{W^q\!\left({\rm Grad};\,\Omega\right)}
				\leq C_4 |\!| g |\!|_{L^q(\Omega)}\,.
		\end{equation}
	Choosing now any constant $C_{5}>C_4$, then for all $g\in L^q(\Omega)$ there exists  
	$v\in W^q\!\left({\rm Grad};\,\Gamma_\tau;\,\Omega\right)$ with ${\rm div\,}v=g$ and 
		\begin{equation}\nonumber
			|\!|v|\!|_{W^q\!\left({\rm Grad};\,\Omega\right)}\leq C_{5}\, |\!| g |\!|_{L^q(\Omega)}\,. 
		\end{equation}
	Thus, Lemma \ref{DivergenceLemma} is completely proved. \hfill $\Box$\\

Now we are able to prove Theorem \ref{Dev-Div-T}, utilizing the idea from Lemma 3.1 in \cite{arnolddouglasgupta:84}.\\

\noindent
{\sl Proof of Theorem \ref{Dev-Div-T}:}
Let $T\in W^q\!\left({\rm Div};\,\Gamma_\nu;\,\Omega\right)$.
	Since by definition  $T= {\rm dev\,}T+ \frac{1}{n}{\rm tr}(T)\cdot {\rm Id}$, 
	it is sufficient to bound  $|\!|  {\rm tr}(T) |\!|_{L^q(\Omega)}$ by the right
	hand side of (\ref{Dev-Div-E}).
 	Employing a corollary of the Hahn-Banach-Theorem, see e.g. \cite[IV.6, Corollary 2]{yoshida:80}, for every $T\in L^q(\Omega)$
	there is a $g\in L^{q'}(\Omega)$ with $|\!| g |\!|_{L^{q'}(\Omega)}=1$ and
		\begin{equation}\nonumber
			|\!| {\rm tr}(T) |\!|_{L^q(\Omega)} = \int_\Omega {\rm tr}(T)\,g\,\dx\,.
		\end{equation}
	Due to Lemma \ref{DivergenceLemma},
	there exists some vector field $v\in W^{q'}\left({\rm Grad};\,\Gamma_\tau;\,\Omega\right)$,
	such that ${\rm div\,}v = g$ and the estimate 
	$|\!| v |\!|_{W^{q'}\left({\rm Grad};\,\Omega\right)}\leq C \,|\!| g |\!|_{L^{q'}(\Omega)}\leq C$
	holds, where $C>0$ does not depend on $g$, $v$ or $T$. Thus, 
		\begin{eqnarray}\label{BoundaryTerms}
				\frac{1}{n}|\!| {\rm tr}(T)|\!|_{L^q(\Omega)} & =&
						\frac{1}{n}\int_\Omega{\rm tr}(T) \,{\rm div\,}v\,\dx
				= \frac{1}{n}\int_\Omega({\rm tr}(T)\cdot{\rm Id}\,, {\rm Grad\,}v)_{\R^{n\times n}}\,\dx  \\ 
				& = &\int_\Omega(T-{\rm dev\,}T, {\rm Grad\,}v)_{\R^{n\times n}}\,\dx
				=  -\int_\Omega {\rm Div\,}T\cdot  v+({\rm dev\,}T\,,{\rm Grad\,}v)_{\R^{n\times n}}\,\dx \nonumber \\
				& \leq & C \left(|\!| {\rm Div\,}T |\!|_{L^q(\Omega)}+|\!| {\rm dev\,}T |\!|_{L^q(\Omega)}\right) \,.\nonumber
		\end{eqnarray}
Note that no boundary terms occur since $v\in W^{q'}\left({\rm Grad};\,\Gamma_\tau;\,\Omega\right)$
and $T\in W^q\!\left({\rm Div};\,\Gamma_\nu;\,\Omega\right)$.
Therefore, (i) and also  (i') are proved. 

Now we prove (ii).
Let $T\in W^q\!\left({\rm Curl};\,\Gamma_\nu;\,\Omega\right)$.
For $n=3$, ${\rm Curl\,}T$ is again a quadratic tensor and the homogeneous tangential trace
is mapped by the Curl operator to the homogeneous normal trace\footnote{
	For the convenience of the reader we give an illustration of this well known fact assuming a
	completely smooth setting: Let $E$ be a row-vector of $T$, $n$ the outward unit normal of $\Omega$, $\times$ 
	the vector product in $\R^3$ and $u$ an 
	arbitrary function with ${\rm supp\,}u\cap\partial\Omega\subset\Gamma_\nu$. 
	Using Gauss Theorem twice  we compute
	\begin{align*}
		\int_{\partial\Omega} (n\cdot {\rm curl\,}E) u\,do &=\int_\Omega {\rm div\,}(u \,{\rm curl\,}E)\dx=
		\int_\Omega {\rm grad\,} u\cdot {\rm curl\,}E\dx = \int_\Omega {\rm grad\,} u\cdot {\rm curl\,}E-({\rm curl\,}{\rm grad\,} u)\cdot E\dx\\
		&=  \int_\Omega {\rm div\,}(E\times {\rm grad\,} u)\dx=\int_{\partial\Omega} n\cdot(E\times{\rm grad\,} u)\,do = 
			\int_{\partial\Omega} (n\times E)\cdot{\rm grad\,} u\,do=0.
	\end{align*}
Since $u$ is arbitrary the normal trace of ${\rm curl\,} E$ is vanishing on $\Gamma_\nu$.
(Using Stokes's theorem the same is proved in one line.)}.
	Thus  ${\rm Curl}\,T$ belongs to $W^q\!\left({\rm Div};\,\Gamma_\nu;\,\Omega\right)$  and furthermore ${\rm Div\,}{\rm Curl\,}T=0$. 
	Now (ii) follows immediately by \eqref{Dev-Div-E} applied to ${\rm Curl\,}T$.
		\hfill$\Box$

	
\section{The DevSym-Curl-inequality}

Sym-Curl-estimates have been established recently in a series of papers by some of the present authors and have been shown to hold
true also for mixed boundary conditions, see \cite{neffpaulywitsch:13a} for $n=3$ 
and \cite{neffpaulywitsch:13b} for arbitrary dimensions.
For these results it is crucial that the domain $\Omega$ allows for the so-called 
Maxwell compactness property, i.e. the compact embedding
of $H\!\left({\rm curl};\,\Gamma_\tau;\,\Omega\right)\cap H\!\left({\rm div};\,\Gamma_\nu;\,\Omega\right)$ into $L^2(\Omega)$,
and the so-called Maxwell approximation property, see \cite{neffpaulywitsch:13b}.
These two properties ensure that the
Helmholtz decomposition (also for tensor fields) holds true, see \cite{neffpaulywitsch:13a,neffpaulywitsch:13b}:
	\begin{equation}\label{HelmholtzDecomposition}
		L^2(\Omega) = {\rm Grad\,}H\!\left({\rm Grad};\,\Gamma_\tau;\,\Omega\right)\oplus\mathcal H(\Omega)\oplus
						{\rm Curl\,}H\!\left({\rm Curl};\,\Gamma_\nu;\,\Omega\right)\,,
	\end{equation}
where $\mathcal H(\Omega)$ is the space of harmonic Dirichlet-Neumann-tensors, i.e., the space of  tensors 
$T$ belonging to 
$H\!\left({\rm Curl};\,\Gamma_\tau;\,\Omega\right)\cap H\!\left({\rm Div};\,\Gamma_\nu;\,\Omega\right)$ 
with ${\rm Curl\,}T=0$ and ${\rm Div\,}T=0$, and $\oplus$ denotes orthogonality in $L^2(\Omega)$.
Due to the Maxwell compactness property, the unit ball in $\mathcal H(\Omega)$ is compact
and hence the space $\mathcal H(\Omega)$ has finite dimension, the dimension depending on topological
properties of the domain.
In consequence of the Maxwell compactness property, a Poincar\'e--type Maxwell estimate is achieved
by a standard indirect argument, i.e.
	\begin{equation}\label{Maxwell-E}
		|\!| T |\!|_{L^2(\Omega)} \leq C_m \left(|\!| {\rm Curl\,}T |\!|_{L^2(\Omega)}+|\!| {\rm Div\,}T |\!|_{L^2(\Omega)}\right)
	\end{equation}
	for all $T\in H\!\left({\rm Curl};\,\Gamma_\tau;\,\Omega\right)\cap H\!\left({\rm Div};\,\Gamma_\nu;\,\Omega\right)$
	perpendicular to $\mathcal H(\Omega)$, see \cite{neffpaulywitsch:13a}.
Both, the Maxwell compactness property and the Maxwell approximation property have 
been proved to be satisfied, if the underlying domain
$\Omega$ has a Lipschitz boundary, and in addition the interface between the two kinds of boundaries
	\begin{equation}\label{InterfaceLipschitz}
			\bar{\Gamma}_\tau\cap\bar{\Gamma}_\nu \quad\text{is also Lipschitz,}
	\end{equation}
see \cite{jochmann:97,jakabmitreamitra:09} and the discussion in \cite{neffpaulywitsch:13a,neffpaulywitsch:13b}.

In order to deal with the influence of possible harmonic Dirichlet-Neumann-tensors, 
in \cite[Definition 10]{neffpaulywitsch:13a} a further  technical condition on the domain $\Omega$ and the topology of $\Omega$
is imposed: 

\begin{definition}\label{DefinitionSliceable}
	$\Omega$ is called sliceable, if there exist a natural number $J\in \N$ and $\Omega_j\subset\Omega$, $j=1,\dots, J$, 
	such that $\Omega\setminus \left(\Omega_1\cup\dots\cup\Omega_J\right)$ is a null set and for $j=1, \dots, J$
		\begin{itemize}
			\item[(i)] 	$\Omega_j$ are open, disjoint and simply connected Lipschitz subdomains of $\Omega$,
			\item[(ii)] 	$\Gamma_{t, j}:={\rm int_{\rm rel}}\left(\bar{\Omega}_{j}\right)\cap \Gamma_\tau\not=\emptyset$,
					if $\Gamma_\tau\not=\emptyset$.
		\end{itemize}
\end{definition}

First we prove:

\begin{lemma}\label{DevSymGrad-T}
	 Let $n\geq 3$ and $\Gamma_\tau\not=\emptyset$
	 or $n=2$ and $\Gamma_\tau=\partial\Omega$.
	Then, there is a constant $C_{D\!S\!G}$, 
	such that for all
	$v\in H\!\left({\rm Grad};\,\Gamma_\tau;\,\Omega\right)$ 
		\begin{equation}\label{DevSymGrad-Estimate}
			|\!| {\rm Grad\,} v|\!|_{L^2(\Omega)}\leq C_{D\!S\!G} |\!| {\rm dev\,}{\rm sym\,} {\rm Grad\,}v|\!|_{L^2(\Omega)}.
		\end{equation}
\end{lemma}

The proof of Lemma \ref{DevSymGrad-T} relies only on the estimate (\ref{KornSecondTraceFree-E}), i.e.,
an improved version of Korn's second inequality, Rellich's selection theorem
and the control of the kernel of ${\rm dev\,}{\rm sym\,}{\rm Grad\,}$ through the boundary condition.
On this account, a representation formula for elements in this kernel is needed,  
which is given in the Appendix of this paper.
The case $n=2$ with full boundary condition is already proved in the Appendix of \cite{neffpaulywitsch:13a}
and a counterexample to (\ref{DevSymGrad-Estimate}) for the case $n=2$ 
without the full boundary condition will be given in Section 6.\\

\noindent	
{\sl Proof:}
	In a first step, we prove
		\begin{equation}\label{KernDevSymGradBoundary}
			\left(v\in H\!\left({\rm Grad};\,\Gamma_\tau;\,\Omega\right)\quad\wedge\quad
			{\rm dev\,}{\rm sym\,}{\rm Grad\,}v=0\right)\quad \Rightarrow\quad v=0\,.
		\end{equation}
	For the case $n=3$ this is already proved e.g. in \cite{Neff_JeongMMS08}. 
	Here we deal with the case of arbitrary space dimension $n\geq3$.
	We utilize the following representation of the kernel which is proved in the Appendix: 
	There are vectors $\bar{v}, \bar{w}\in \R^n$, 
	a real number $\bar{\p}\in \R$ and a skew-symmetric matrix $\bar{A}\in \mathfrak{so}(n)$, such that
	\begin{eqnarray}\label{RepresentationFormulaFunction_T}
		v(x) 	& = &  \p(x)\,x 
		-\frac{1}{2}|x|^2\bar{w}+\bar{A}x+\bar{v} \: ,\\ 	\label{RepresentationFormulaGradient_T}
			{\rm grad\,} v(x) 	& = &	\p(x)\,{\rm Id}+A(x)\,,
	\end{eqnarray}
holds for all $x\in\bar{\Omega}$, where 
		\begin{equation}\label{RepresentationpA_T}
			\p(x) 		 =  \bar{w}\cdot x +\bar{\p} \: ,\quad
			A_{ij}(x)		=  \sum_{k=1}^n\bar{a}_{ijk}x_k+\bar{A}_{ij}
		\end{equation}
	and
		\begin{equation}\label{aijk}
				\bar{a}_{ijk} = 
					\begin{cases}
						0 & \text{if}\quad i\not=j, i\not= k, k\not=j, \\
						0 & \text{if}\quad i=j, \\
						\bar{w}_j  & \text{if}\quad k=i,  i\not=j, \\
						-\bar{w}_i & \text{if}\quad k=j,  k\not=i\,.
					\end{cases}
		\end{equation}
	In particular, $A(x)$ is skew-symmetric and the dimension of the kernel of ${\rm dev\,}{\rm sym\;}{\rm Grad\,}$ is 
	$(n+1)(n+2)/2$.
	Due to this formula, $v$ is a smooth vector field on $\bar{\Omega}$.
	Let $x\in\Gamma_\tau$ and $\tau\in\R^n$, $\tau\not=0$, tangential to $\Gamma_\tau$ in $x$. 
	Since $v\in H\left({\rm Grad};\,\Gamma_{\tau};\,\Omega\right)$, we have
	${\rm grad\,} v\in H\left({\rm Curl};\,\Gamma_{\tau};\,\Omega\right)$, i.e.
	$$ {\rm grad\,} v(x)\,\tau=0\;.$$
	Therefore, if $x\in\Gamma_{\tau}$, then ${\rm grad\,} v(x)$ does not have full rank. 
	By (\ref{RepresentationFormulaGradient_T}) and since $\tau\cdot A(x)\,\tau=0$
		\begin{equation}\label{Betrag}
			0=| {\rm grad\,} v(x) \,\tau |^2= \p^2(x) |\tau|^2+|A(x) \tau|^2
		\end{equation}
	holds with $\p$ and $A$ from (\ref{RepresentationpA_T}). 
	Hence, $\p(x)=0$.
Therefore,	(\ref{RepresentationpA_T})	 implies necessarily 
		\begin{equation}\label{Hypersurface}
			0 = \p(x)=\bar{w}\cdot x +\bar{\p}\quad\text{for all}\quad x\in\Gamma_\tau\,.
		\end{equation}
	On the other hand, if $\p(x)=0$, then ${\rm grad\,} v(x)$ has not full rank, since $A(x)$ is skew-symmetric.
	Thus, for all $x\in\Gamma_\tau$ the matrix ${\rm grad\,} v(x)$  does not have full rank, if and only if (\ref{Hypersurface}) holds.	
	If $\bar{w}\neq0$, then by \eqref{Hypersurface} $\Gamma_\tau\subset E$, where $E$ denotes
	the affine hypersurface defined by equation (\ref{Hypersurface}). On the other hand, for all
	$x\in\Gamma_\tau\subset E$, due to the representation formula (\ref{RepresentationFormulaFunction_T})
	and (\ref{Hypersurface}), we get
		\begin{equation}\label{QuadraticEquation}
			v(x)=	-\frac{1}{2}|x|^2\bar{w}+\bar{A}x+\bar{v}=0\,,
		\end{equation}
	 describing for $\bar{w}\not=0$ a quadratic surface and not a hypersurface. 
	 This proves $\bar{w}=0$ and hence $\p=\bar{\p}=0$.
	Consequently, on $\Gamma_\tau$ 
		\begin{equation}\label{LinearEquation}
			v(x)=\bar{A}x+\bar{v}=0\,,
		\end{equation}
	yielding $\bar{A}=0$ and $\bar{v}=0$, since otherwise 
	the solution set of  (\ref{LinearEquation}) is an affine surface with  co-dimension ${\rm codim}\geq 2$, 
	recall that $\bar{A}$ is skew-symmetric.
	But such a surface cannot contain an open and non-empty subset of a Lipschitz-continuous boundary. 
	Therefore (\ref{KernDevSymGradBoundary}) is proved.

	In the second step we utilize \ref{KornSecondTraceFree-E} from \cite[Theorem 1.1]{dain:06} or \cite{Neff_JeongMMS08}
	and carry out the usual conclusion by contradiction.
	Assume the estimate (\ref{DevSymGrad-Estimate}) is false, then there exists a sequence 
	$(v_j)\subset H\!\left({\rm Grad};\,\Gamma_\tau;\,\Omega\right)$ with
	$|\!| {\rm Grad\,} v_j |\!|_{L^2(\Omega)}=1$ and
		\begin{equation}\label{devsymGradtozero}
			|\!| {\rm dev\,}{\rm sym\;}{\rm Grad\,} v_j |\!|_{L^2(\Omega)}<\frac{1}{j}
		\end{equation}
	for all $j\in\N$. According to  (\ref{KornSecondTraceFree-E})  the sequence of norms $|\!| v_j |\!|_{L^2(\Omega)}$ is 
	bounded from below, i.e., there exists $J\in\N$ and a constant $C>0$, such that
		\begin{equation}\label{boundvj}
			|\!| v_j |\!|_{L^2(\Omega)}\geq C	\quad\text{for all}\quad j\geq J.
		\end{equation}
	Utilizing Poincar\'{e}'s inequality and $|\!| {\rm Grad\,} v_j |\!|_{L^2(\Omega)}=1$,
	the sequence $(v_j)$ is bounded in $H\!\left({\rm Grad};\,\Omega\right)$. Employing Rellich's selection theorem
	there is a subsequence of $(v_j)$, again called
	$(v_j)$, and  $v\in H\!\left({\rm Grad};\,\Gamma_\tau;\,\Omega\right)$ such that 
		\begin{eqnarray}\label{vjstronglytov}
			v_j & \rightarrow & v \quad\text{strongly in}\quad L^2(\Omega)\,, \\ \nonumber
			{\rm Grad\,} v_j& \rightharpoonup & {\rm Grad\,}v\quad\text{weakly in}\quad L^2(\Omega)\,.
		\end{eqnarray}
	Hence, ${\rm dev\,}{\rm sym\,}{\rm Grad\,}v_j$ converges weakly to ${\rm dev\,}{\rm sym\,}{\rm Grad\,}v$
	and due to weak lower semi-continuity of the norm and (\ref{devsymGradtozero}) we conclude
		\begin{equation}\nonumber
			|\!| {\rm dev\,}{\rm sym\,}{\rm Grad\,} v |\!|_{L^2(\Omega)}
			\leq\liminf_{j\to\infty}|\!| {\rm dev\,}{\rm sym\,}{\rm Grad\,} v_j |\!|_{L^2(\Omega)}=0\,.
		\end{equation}
	According to (\ref{KernDevSymGradBoundary}), this implies $v=0$, in contradiction to (\ref{boundvj}) und (\ref{vjstronglytov}).
	Therefore, Lemma \ref{DevSymGrad-T} is proved.\hfill$\Box$\\

Now, we can prove the DevSym-Curl-inequality:

\begin{theorem}\label{DevSym-Curl-T}
	Let $n\geq 3$, $\Omega\subset \R^n$ be a slicable domain and $\Gamma_{\tau}\neq\emptyset$.
	Then, there is a positive constant $C_{D\!S\!C}$, 
	such that  for all
	$T\in H\!\left({\rm Curl};\,\Gamma_\tau;\,\Omega\right)$
		\begin{equation}\label{DevSym-Curl-E}
				|\!| T |\!|_{L^2(\Omega)}	 \leq 	 C_{D\!S\!C} \left(|\!| {\rm dev\,}{\rm sym\,} T|\!|_{L^2(\Omega)}
													+|\!| {\rm Curl\,} T|\!|_{L^2(\Omega)}\right).
     \end{equation}
\end{theorem}

We note that Theorem \ref{DevSym-Curl-T} remains true if $n=2$ and $\Gamma_{\tau}=\partial\Omega$
since Lemma \ref{DevSymGrad-T} holds in this case as well. Moreover, with \eqref{DevSym-Curl-E} also
$$	|\!| T |\!|_{L^2(\Omega)}+|\!|  {\rm Curl\,}T |\!|_{L^2(\Omega)}
			\leq C_{D\!S\!C}\left(|\!| {\rm dev\,}{\rm sym\,} T|\!|_{L^2(\Omega)}
			+|\!| {\rm Curl\,} T|\!|_{L^2(\Omega)}\right)$$
holds.\\

\noindent
{\sl Proof:}
	We combine the proof of the Sym-Curl-inequality  (\ref{Sym-Curl-E}) 
	from the papers \cite{neffpaulywitsch:13a,neffpaulywitsch:13b} 
	with Lemma \ref{DevSymGrad-T}. 
	Let $T\in H\!\left({\rm Curl};\,\Gamma_\tau;\,\Omega\right)$. Using the Helmholtz decomposition from 
	\cite{neffpaulywitsch:13b} we have the orthogonal sum 
		\begin{equation}\nonumber
			T = R+S\in H\!\left({\rm Curl_0};\,\Gamma_\tau;\,\Omega\right)
			\oplus {\rm Curl\,}H\!\left({\rm Curl};\,\Gamma_\nu;\,\Omega\right)\,,
		\end{equation} 
	where $ R \in\left({\rm Curl_0};\,\Gamma_\tau;\,\Omega\right)$, if and only if
	$R\in H\!\left({\rm Curl};\,\Gamma_\tau;\,\Omega\right)$ and ${\rm Curl\,}R=0$.
	Note, that in general $R\in H\!\left({\rm Curl_0};\,\Gamma_\tau;\,\Omega\right)$
	does not imply $R={\rm Grad\,}v$ with
	$v\in H\!\left({\rm Grad};\,\Gamma_\tau;\,\Omega\right)$, since, depending on topological
	properties of the domain $\Omega$, some harmonic-Dirichlet-Neumann tensor fields could be involved.
	In order to deal with this possibility, we slice the domain $\Omega$ according to Definition \ref{DefinitionSliceable}
	and set
		\begin{equation}\nonumber
				R=\sum_{j=1}^J \chi_j R_j \,,
		\end{equation} 
	where $R_j:= R|_{\Omega_j}$ and $\chi_j$ is the indicator-function of $\Omega_j$.
	In the proofs of \cite[Lemmas 9 and 12]{neffpaulywitsch:13a} it is shown, that there are
	non-empty and relatively open connected
	subsets $\widetilde{\Gamma}_{\tau,j}\subset \Gamma_{\tau,j}$ and vector fields
	$v_j\in H\!({\rm Grad};\,\widetilde{\Gamma}_{\tau,j};\,\Omega_j)$
	such that ${\rm Grad\,}v_j= R_j$. Now we apply (\ref{DevSymGrad-Estimate}) to $v_j$ and get
	\begin{eqnarray}\label{EstimateR}
		|\!| T |\!|^2_{L^2(\Omega)}	& =  &|\!| R |\!|^2_{L^2(\Omega)}+|\!| S |\!|_{L^2(\Omega)}^2
		= \sum_{j=1}^J	|\!| R_j |\!|_{L^2(\Omega_j)}^2 + |\!| S |\!|_{L^2(\Omega)}^2
		\\ \nonumber
		& \leq & 
		  C \sum_{j=1}^J	|\!| {\rm dev\,}{\rm sym\,}R_j|\!|_{L^2(\Omega_j)}^2
												+|\!| S |\!|_{L^2(\Omega)}^2 
		 =  C\,|\!| {\rm dev\,}{\rm sym\,}R|\!|_{L^2(\Omega)}^2+|\!| S |\!|_{L^2(\Omega)}^2 \nonumber \\
		& \leq & C\,|\!| {\rm dev\,}{\rm sym\,}T|\!|_{L^2(\Omega)}^2
		+ C\,|\!| {\rm dev\,}{\rm sym\,}S|\!|_{L^2(\Omega)}^2
		+|\!| S |\!|_{L^2(\Omega)}^2 \nonumber \\ 
		& \leq & C\,|\!| {\rm dev\,}{\rm sym\,}T|\!|_{L^2(\Omega)}^2
		+ C \, |\!| S |\!|_{L^2(\Omega)}^2 \,.\nonumber
	\end{eqnarray}
	Concerning the $S$-part, we note that ${\rm Curl\,}T={\rm Curl\,}S$ and 
	$S\in H\left({\rm Curl};\,\Gamma_{\tau};\,\Omega\right)$ 
	since $T$ and $R$ belong to $H\left({\rm Curl};\,\Gamma_{\tau};\,\Omega\right)$.
	Moreover, since 
	$${\rm Curl}\,H\left({\rm Curl};\,\Gamma_{\nu};\,\Omega\right)\subset
	H\left({\rm Div}_{0};\,\Gamma_{\nu};\,\Omega\right)\cap\mathcal H(\Omega)^{\bot}$$
	we even have $S\in H\left({\rm Curl};\,\Gamma_{\tau};\,\Omega\right)\cap
	H\left({\rm Div}_{0};\,\Gamma_{\nu};\,\Omega\right)\cap\mathcal H(\Omega)^{\bot}$.
	By means of  the Maxwell inequality (\ref{Maxwell-E}) and since ${\rm Div}\,S=0$ we estimate
		\begin{eqnarray}\label{S-E}
			|\!| S |\!|_{L^2(\Omega)}	
									& \leq & 	C_m\, |\!| {\rm Curl\,}S |\!|_{L^2(\Omega)}
									 =  	C_m\, |\!| {\rm Curl\,}T |\!|_{L^2(\Omega)} \,.
		\end{eqnarray}
	Combining (\ref{S-E}) and (\ref{EstimateR}) yields
		\begin{eqnarray}\nonumber
			|\!| T |\!|_{L^2(\Omega)} & \leq & C\,|\!| {\rm dev\,}{\rm sym\,}T |\!|_{L^2(\Omega)}
												+ C\, |\!| {\rm Curl\,}T |\!|_{L^2(\Omega)}\,,
			\end{eqnarray}
completing the proof.
\hfill$\Box$

\section{The Sym-DevCurl- and DevSym-DevCurl-inequalities}

Now, we combine the Dev-Div-inequality with the Sym-Curl-inequality
and the DevSym-Curl-inequality. For this, we need $n=3$ since only then
${\rm Curl\,} T$ is again quadratic.

\begin{theorem}\label{DevSym-DevCurl-T}
	Let $\Omega\subset \R^3$ be a slicable domain and $\Gamma_{\tau}\neq\emptyset$.
	Then, there are positive constants $C_{S\!D\!C}$ and $C_{D\!S\!D\!C}$, 
	such that  for all
	$T\in H\!\left({\rm Curl};\,\Gamma_\tau;\,\Omega\right)$
		\begin{equation}\label{DevSym-DevCurl-E}
			\begin{array}{rclc}
				|\!| T |\!|_{L^2(\Omega)}	& \leq &	 C_{S\!D\!C} \left(|\!| {\rm sym\,} T|\!|_{L^2(\Omega)}
													+|\!| {\rm dev\,}{\rm Curl\,} T|\!|_{L^2(\Omega)}\right),\\
				|\!| T |\!|_{L^2(\Omega)}	& \leq &	 C_{D\!S\!D\!C} \left(|\!| {\rm dev\,}{\rm sym\,} T|\!|_{L^2(\Omega)}
													+|\!| {\rm dev\,}{\rm Curl\,} T|\!|_{L^2(\Omega)}\right)
			\end{array}
		\end{equation}
		and
		\begin{equation}\label{DevSym-DevCurl-E+Curl}
			\begin{array}{rclc}
				|\!| T |\!|_{L^2(\Omega)} + |\!| {\rm Curl\,}T |\!|_{L^2(\Omega)}	
				& \leq &	 C_{S\!D\!C} \left(|\!| {\rm sym\,} T|\!|_{L^2(\Omega)}
													+|\!| {\rm dev\,}{\rm Curl\,} T|\!|_{L^2(\Omega)}\right),\\
				|\!| T |\!|_{L^2(\Omega)} + |\!| {\rm Curl\,}T |\!|_{L^2(\Omega)}	
				& \leq &	 C_{D\!S\!D\!C} \left(|\!| {\rm dev\,}{\rm sym\,} T|\!|_{L^2(\Omega)}
													+|\!| {\rm dev\,}{\rm Curl\,} T|\!|_{L^2(\Omega)}\right).
			\end{array}
		\end{equation}
\end{theorem}

\noindent
{\sl Proof:}
Combine Theorem \ref{DevSym-Curl-T} with Theorem \ref{Dev-Div-T} (ii).
\hfill$\Box$

\section{Kernels and counterexamples}

It is illuminating to see, how the kernels of the inequalities are controlled on, say, the space of smooth
compactly supported tensor fields. Of course, some of the given arguments are well known.
In the following we always assume that $T$ is such a smooth tensor field with compact support
in $\Omega\subset \R^3$. 

\subsection{The kernel of the Dev-Div-inequality}

Consider some $T$ in the kernel of the Dev-Div-inequality, i.e., ${\rm dev\,}T=0$ and ${\rm Div\,} T=0.$ 
Since ${\rm dev\,}T=0$ we have $T=\p\cdot {\rm Id}$.
But therefore ${\rm Div\,}T={\rm grad\,}\p=0$ and we conclude $\p=const$. Since $\p$ and $T$ are
compactly supported, $\p=0$ and $T=0$ in $\Omega$.

\subsection{The kernel of the Sym-Curl-inequality}

Consider some $T$ in the kernel of the Sym-Curl-inequality, i.e., ${\rm sym\;}T=0$ and ${\rm Curl\,}T=0$. 
Since ${\rm sym\;}T=0$ we conclude $T(x)=A(x)\in\mathfrak{so}(3)$, say
	\begin{equation}\label{SkewA}
		A=\left(\begin{array}{ccc}
				0     & -a_3 & a_2 \\
				a_3  & 0	 & -a_1 \\
				-a_2 & a_1 & 0
			\end{array}\right)\quad\text{and}\quad
	\end{equation}
	\begin{equation}\label{CurlSkewA} 
		{\rm Curl\,}A=\left(\begin{array}{ccc}
				\partial_2 a_2+\partial_3 a_3     & -\partial_1a_2 & -\partial_1a_3 \\
				-\partial_2a_1  &  \partial_3a_3+\partial_1a_1 & -\partial_2a_3 \\
				-\partial_3a_1 & -\partial_3a_2 & \partial_1a_1+\partial_2a_2
			\end{array}\right)
	\end{equation}
with a smooth and compactly supported vector field $a=(a_1, a_2, a_3)^\top $.
Hence ${\rm Curl\,}A=0$ implies  ${\rm Grad\,}a=0$ and thus $a=0$ and $T=A=0$, see also \cite{MuenchNeff:08}.

\subsection{The kernel of the DevSym-DevCurl-inequality}

Regarding the DevSym-DevCurl-inequality the situation gets more involved.
Let us assume ${\rm dev\,}{\rm sym\;}T=0$ and ${\rm dev\,}{\rm Curl\,}T=0$.
Then
\begin{align}
T(x)&=\p(x)\cdot {\rm Id}+A(x),\label{KernelDevSym}\\
{\rm Curl\,}(\p(x)\cdot{\rm Id})+{\rm Curl\,}A(x)={\rm Curl\,}T(x)&=\q(x)\cdot {\rm Id}\,\label{Kerneldevsym-devcurl}
\end{align}
with smooth and compactly supported functions $\p,\q$ 
and with $a$, $A$ as above. 
Now 
	\begin{equation}\label{Curlp}
		{\rm Curl\,}(\p\cdot{\rm Id})=
			\left(\begin{array}{ccc}
					0 & \partial_3\p & -\partial_2\p \\
					-\partial_3 \p & 0 & \partial_1 \p \\
					\partial_2\p & -\partial_1 \p & 0 
				\end{array}\right)
	\end{equation}
	is a skew-symmetric matrix. Therefore, ${\rm sym\;}{\rm Curl\,}A=\q\cdot {\rm Id}$
	and hence by \eqref{CurlSkewA} 
		\begin{equation}\label{Anticommutator}
				\partial_1a_2+\partial_2a_1 = 
				\partial_2a_3+\partial_3a_2 = 
				\partial_3a_1+\partial_1a_3 = 0
		\end{equation}
	and
		\begin{equation}
				\partial_2a_2+\partial_3a_3 =\partial_3a_3+\partial_1a_1 =\partial_1a_1+\partial_2a_2 =\q\,.
		\end{equation}
	The second series of equations yields 
	\begin{equation}\label{Diagonalelemente}
			\partial_1a_1 =\partial_2a_2 =\partial_3a_3 =\frac{\q}{2}\qquad\text{as well as}\qquad
				2\,{\rm div\,}a=3 \,\q\,.
	\end{equation}
	 By means of comparison of the skew-symmetric parts of equation (\ref{Kerneldevsym-devcurl}),
	utilizing (\ref{SkewA}) and (\ref{Curlp}),  we
	conclude that
		\begin{equation}\label{gradientp}
				{\rm grad\,}\p = \left(\begin{array}{c} \partial_2a_3 \\ \partial_3a_1 \\  \partial_1a_2
		\end{array} \right)
	\end{equation}
and thus, employing (\ref{Anticommutator})
		\begin{equation}\label{curlgrad}
			0 = {\rm curl\,}{\rm grad\,}\p= 
			{\rm curl\,}\left(\begin{array}{c} \partial_2a_3 \\ \partial_3a_1 \\  
		\partial_1a_2\end{array} \right)=
			-\left(\begin{array}{c} (\partial_2^2+\partial_3^2)a_1 \\
								(\partial_3^2+\partial_1^2)a_2 \\
								(\partial_1^2+\partial_2^2)a_3\end{array} \right)\,.
		\end{equation}
	With (\ref{curlgrad}) and (\ref{Diagonalelemente}) we obtain
	\begin{equation}\label{Laplace}
			\Delta a = -\frac{1}{2}{\rm grad\,}\q=-\frac{1}{3}{\rm grad\,}{\rm div\,}a\,.
	\end{equation}
	Furthermore, due to	(\ref{Anticommutator}) 
		\begin{equation}\nonumber
			{\rm curl\,}a = 2\left(\begin{array}{c} \partial_2a_3 \\ \partial_3a_1 \\  \partial_1a_2\end{array} \right)
		\end{equation}		
		and employing (\ref{curlgrad}) it follows that ${\rm curl\,}{\rm curl\,}a=0$.
		The combination of this fact with (\ref{Laplace}) and the identity 
		${\rm grad\,}{\rm div\,}-{\rm curl\,}{\rm curl\,}=\Delta$ yields
			\begin{equation}\nonumber
				{\rm grad\,}{\rm div\,}a
				=\Delta a
				= -\frac{1}{3}{\rm grad\,}{\rm div\,}a
			\end{equation}
		and thus ${\rm grad\,}{\rm div\,}a=\Delta a=0$.	Since this Poisson equation is uniquely
		solvable we conclude $a=0$ and $A=0$,  and utilizing (\ref{gradientp}) also $\p=const$. 
		Hence, $\p=0$ yielding $T=0$.  

\subsection{There are no DevSym-DevSymCurl- or DevSym-SymCurl-inequalities}

	Choose $\p\in C^\infty_0(\Omega;\,\R)$ and set $T:=\p\cdot{\rm Id}$. 
	Then ${\rm dev\,}{\rm sym\;}T=0$ and, according
	to  (\ref{Curlp}), ${\rm sym\;}{\rm Curl\,}T=0$. Therefore, such inequalities have to be false.
	
\subsection{There is no Sym-Div-inequality} 
	
	Choose $\p\in C^\infty_0(\Omega;\,\R)$ and set $a:={\rm grad\,}\p$ and
	define $A$ according to (\ref{SkewA}). 
	Then we have ${\rm Div\,}A=-{\rm curl\,}a = -{\rm curl\,}{\rm grad\,}\varphi=0$ 
	and ${\rm sym\;}A=0$. Therefore, such an inequality is false.

\subsection{The DevSymGrad-inequality is false for $n=2$}

	As already announced in the introduction, now we show that in the case $n=2$ 
	the trace-free version of Korn's first inequality with only
	partial boundary condition is false. 
	This is remarkable, since the kernel of the inequality is already controlled by a partial boundary condition.
	In fact, if a function is in the kernel, then it is holomorphic in $\Omega$. But if a holomorphic function 
	vanishes on some part of the boundary it has to vanish on the whole of $\Omega$. This shows that having a norm 
	on the space under consideration is only necessary for the validity of an inequality.
	The construction of our counterexample is taken from
	\cite{pompe:10} and in that paper it served as a counterexample to a version of Korn's first inequality, see \cite{neff00b},
	with non-constant (rotation) coefficients, see also \cite{neff_pompe13}. 
	For the convenience of the reader we introduce this example in detail,
	thereby we exactly follow \cite{pompe:10}.

		We identify $\R^2$ with $\C$ via standard notation $z=x+iy$. We also use the standard notation for the polar coordinates
		$(x, y)=r\,(\cos t,\sin t)$.
		 Consider the sequence 
		 $$u_n(x, y)=xz^n$$ 
		 on the half disk $\Omega = \{z \,:\,|z|<1,\,x>0\}$.
		As $\Gamma_{\tau}$ we choose $\{z\in\partial\Omega \,:\,x=0\}$. 
		Then, of course, each of the mappings $u_n$ vanishes
		on $\Gamma_\tau$. We first compute ${\rm grad\,} u_n(x, y).$ Since 
		$$(z^n)' = nz^{n-1} = nr^{n-1}(\cos(nt-t)+i\sin(nt-t)),$$ 
		we obtain
			\begin{equation}\nonumber
				{\rm grad\,} z^n= nr^{n-1}\left(\begin{array}{cc}
											\cos(nt-t) & -\sin(nt-t) \\
											\sin(nt-t)  & \cos(nt-t)	
										\end{array}\right)\,.
			\end{equation}
		Therefore, we have
		\begin{eqnarray}\label{AbleitungUn}
			{\rm grad\,} u_n(x, y) & = & (x{\rm grad\,}_x(z^n)+z^n, x{\rm grad\,}_y(z^n))\\ \nonumber
					    &	 =  & r^{n}\left(\begin{array}{cc}
											\cos(nt) & 0 \\
											\sin(nt)  & 0	
										\end{array}\right)+
								nr^{n}\cos t \left(\begin{array}{cc}
											\cos(nt-t) & -\sin(nt-t) \\
											\sin(nt-t)  &  \cos(nt-t)	
										\end{array}\right)
		\end{eqnarray}
and hence
		\begin{eqnarray}\nonumber
			|{\rm grad\,} u_n|^2 	& = & r^{2n} + 2n^2r^{2n}\cos^2 t + 2r^{2n}n\,\cos t \,
								(\cos(nt)\,\cos(nt-t)+\sin(nt)\,\sin(nt-t))\\ \label{GradUn}
					 	& = & r^{2n} + 2r^{2n}(n^2+n)\cos^2t.
		\end{eqnarray}	
	Taking into account that
		\begin{equation}\nonumber
			\int_{-\pi/2}^{\pi/2}\cos^2 t \,dt = \frac{\pi}{2}\,,
		\end{equation}
	we obtain
	\begin{eqnarray*}
		\int_\Omega |{\rm grad\,} u_n|^2 & = & \int_0^1\int_{-\pi/2}^{\pi/2}
									r(r^{2n} + 2r^{2n}(n^2+n)\cos^2t)\,dt\,dr\\
								& = & \pi(n^2+n+1)\int_0^1r^{2n+1}\,dr=\pi \frac{n^2+n+1}{2n+2}\,.	
	\end{eqnarray*}
	Now, we use this construction as a counterexample for the DevSymGrad-inequality:
	Switching back to our notation we have
	\begin{equation}\label{LimGradUn}
		\lim_{n\to\infty} |\!|{\rm Grad\,}u_n|\!|_{L^2(\Omega)} = \infty\,.
	\end{equation}
	On the other hand, inspection of formula (\ref{AbleitungUn}) yields
		\begin{eqnarray*}
			{\rm sym\;}{\rm Grad\,}u_n & = & r^n \left(\begin{array}{cc}
									\cos(nt) & \frac{1}{2}\sin(nt)\\
									\frac{1}{2}\sin(nt)  & 0
									\end{array}\right)
									+n r^n\cos t \left(\begin{array}{cc}
									\cos(nt-t) & 0 \\
									0  & \cos(nt-t)
									\end{array}\right) \,, \\
			{\rm dev}_2 \, {\rm sym\;}{\rm Grad\,}u_n  & = &  
				\frac{1}{2}r^n \left(\begin{array}{cc}
									\cos(nt) & \sin(nt)\\
									\sin(nt)  & -\cos(nt)
									\end{array}\right) \, ,
		\end{eqnarray*}
		where ${\rm dev}_2 \, X = X - \frac{1}{2} \tr X\cdot {\rm Id}$ denotes the two-dimensional deviator.
		Hence,
		\begin{equation}\nonumber
			|{\rm dev}_2 \, {\rm sym\;}{\rm Grad\,}u_n|^2 = \frac{1}{2} r^{2n}
		\end{equation}
		and thus 
			\begin{equation}\nonumber
			|\!| {\rm dev}_2 \, {\rm sym\;}{\rm Grad\,}u_n |\!|^2_{L^2(\Omega)}=\frac{\pi}{4n+4}\,,
			\end{equation}
		converging to zero in the limit $n\to\infty$ in contradiction to (\ref{LimGradUn}) and (\ref{DevSymGrad-E}).

The fact that the DevSymGrad inequality does not hold in the two-dimensional case is due to the
special form of the ${\rm dev}$ operator in this case. If we instead view the plane symmetric gradient
as an object in three dimensions and apply the standard ${\rm dev}$ operator for $n=3$ (simply
denoted by ${\rm dev}$ in the sequel), then we obtain
	\begin{eqnarray*}
			{\rm dev}\,{\rm sym\;}{\rm Grad\,}u_n  & = &  
				r^n \left(\begin{array}{ccc}
				\frac{2}{3}\cos(nt) & \frac{1}{2}\sin(nt) & 0 \\
				\frac{1}{2}\sin(nt)  & -\frac{1}{3}\cos(nt) & 0 \\
				0 & 0 & -\frac{1}{3}\cos(nt)
				\end{array}\right) \\
			& + & n r^n \cos t \cos(nt-t)\left(\begin{array}{ccc}
				\frac{1}{3}& 0 & 0 \\
				0  & \frac{1}{3} & 0 \\
				0 & 0 & -\frac{2}{3}
				\end{array}\right) \: .
	\end{eqnarray*}
This implies
		\begin{eqnarray*}
			|{\rm dev\,}{\rm sym\;}{\rm Grad\,}u_n|^2 &=& r^{2n}
			\left( \frac{2}{3} \cos^2 (nt) + \frac{1}{2} \sin^2 (nt) \right)
			+ \frac{2}{3} n^2 r^{2n} \cos^2 t \cos^2 (nt-t) \\
			&\geq& \frac{2}{3} n^2 r^{2n} \cos^2 t \cos^2 (nt-t)
		\end{eqnarray*}
		and, for $n > 2$,
			\begin{equation}\nonumber
				|\!| {\rm dev^{3D}\,}{\rm sym\;}{\rm Grad\,}u_n |\!|^2_{L^2(\Omega)} \geq
				\frac{2}{3} n^2 \frac{\pi}{4} \frac{1}{2n+2} = \frac{\pi n^2}{12(n+1)}\,,
			\end{equation}
where we used the fact that
		\begin{equation}\nonumber
			\int_{-\pi/2}^{\pi/2} \cos^2 t \cos^2 (nt-t) \,dt = \frac{\pi}{4}
		\end{equation}
holds for $n > 2$. This means that
$|\!| {\rm dev\,}{\rm sym\;}{\rm Grad\,}u_n |\!|^2_{L^2(\Omega)} \rightarrow \infty$ for
$n \rightarrow \infty$ in concordance with Lemma \ref{DevSymGrad-T}

\section{Applications}

In this section we will present some prototype applications where the new inequalities may be used to establish coercivity of the models. 

\subsection{Infinitesimal incompressible elasticity}

Historically, inequalities like the one in Theorem \ref{Dev-Div-T} first appeared in the context of mixed
stress-displacement formulations of linear elasticity in the incompressible limit (cf.
\cite{arnolddouglasgupta:84}). The result in \cite{arnolddouglasgupta:84} is stated in two dimensions
assuming vanishing average trace (see also \cite[Sect. VII.2]{brezzifortin:91}. It is generalized in
\cite{carstensendolzmann:98} using a different
argument assuming only that the identity tensor is eliminated by some constraint. In the
incompressible limit, the mixed variational formulation of linear elasticity turns into the problem of
finding some $\sigma \in H (\Div ; \Gamma_\nu ; \Omega)$, $u \in L^2 (\Omega;\R^n)$ and
$\gamma \in L^2 (\Omega;\so(n))$, such that
\begin{equation}
  	\begin{split}
    		({\rm dev} \: \sigma , \tau)_{L^2 (\Omega)} + ( u , \Div \tau )_{L^2 (\Omega)}
		+ ( \gamma , {\rm skew} \: \tau )_{L^2 (\Omega)} & = 0 \: , \\
	  	(\Div \sigma , v)_{L^2 (\Omega)} + (f , v )_{L^2 (\Omega)} & = 0 \: , \\
		( {\rm skew} \: \sigma , \eta )_{L^2 (\Omega)} & = 0
	\end{split}
\end{equation}
holds for all $\tau \in H (\Div ; \Gamma_\nu ; \Omega)$, $v \in L^2 (\Omega;\R^n)$ and
$\eta \in L^2 (\Omega;\so(n))$. This saddle-point problem may be viewed as the
Karush-Kuhn-Tucker system associated with minimizing the elastic energy with respect to the
stresses subject to momentum balance and symmetry as constraints. Its well-posedness relies on the
estimate in Theorem \ref{Dev-Div-T}. The same is true for the stress-displacement first-order system
least squares approach studied in \cite{caistarke:04}.

\subsection{Pseudostress formulation of stationary Stokes equations}

Here, the following formulation of the stationary Stokes equations is considered:
For some given $f : \Omega \rightarrow \R^3$ find the pressure $p : \Omega \rightarrow \R$,  the velocity
$u : \Omega \rightarrow \R^3$ and the stress $\sigma : \Omega \rightarrow \R^{3 \times 3}$
such that the first-order system
$$\sigma-\mu\sym{\rm grad\,} u+p \, \id=0,\qquad\Div\sigma=f,\qquad {\rm div} \, u=0$$
holds in $\Omega$. This system is obviously equivalent to
$$\dev\sigma-\mu\sym{\rm grad\,} u=0,\qquad\Div\sigma=f,$$
where the pressure $p$ has been eliminated and can be computed afterwards as
$p=-\tr\sigma/3$. For this first-order system, a
least squares formulation based on minimizing the quadratic functional
\begin{equation}
  \|\dev\sigma-\mu\sym{\rm grad\,} u\|_{L^2(\Omega)}^2 + \|\Div\sigma-f\|_{L^2(\Omega)}^2
\end{equation}
with respect to $u$ and $\sigma$ may be used. In order to obtain a coercivity result for this functional, let us first investigate the mixed terms arising in the first part of the functional, leading to
\begin{align*}
  ( \dev\sigma,\sym{\rm grad\,} u )_{L^2 (\Omega)}
  & = ( \sym\dev\sigma,{\rm grad\,} u )_{L^2 (\Omega)}
    = ( \sym\sigma - \frac{1}{3} \tr\sigma \id,{\rm grad\,} u )_{L^2 (\Omega)} \\
  & = ( \sigma , {\rm grad\,} u )_{L^2 (\Omega)}
    - ( \skew \sigma , {\rm grad\,} u )_{L^2 (\Omega)}
    - \frac{1}{3} ( \tr \sigma , {\rm div} \: u )_{L^2 (\Omega)} \\
  & = - ( \Div \sigma , u )_{L^2 (\Omega)}
    - ( \skew \sigma , {\rm grad\,} u )_{L^2 (\Omega)}
    - \frac{1}{3} ( \tr \sigma , {\rm div} \: u )_{L^2 (\Omega)} \: ,
\end{align*}
if we assume proper boundary conditions on $\sigma$ and $u$,
justifying the partial integration without boundary terms, i.e.,
$\sigma \in H (\Div ; \Gamma_\nu ; \Omega)$ and  $u \in H (\Grad ; \Gamma_\tau ; \Omega)$.
This implies, for arbitrary $\delta \in (0,1)$,
\begin{equation}
	\begin{split}
		2 \mu ( \dev\sigma,\sym{\rm grad\,} u )_{L^2 (\Omega)}
     		& \leq \delta \left( \mu^2 \| u \|_{L^2 (\Omega)}^2 + \mu^2 \| {\rm grad\,} u \|_{L^2 (\Omega)}^2
                                                    + \frac{1}{3} \| \tr \sigma \|_{L^2 (\Omega)}^2 \right) \\
                   & \qquad + \frac{1}{\delta} \left( \| \Div \sigma \|_{L^2 (\Omega)}^2
                                                    + \| \skew \sigma \|_{L^2 (\Omega)}^2
                                                    + \frac{1}{3} \mu^2 \| {\rm div} \: u \|_{L^2 (\Omega)}^2 \right) \: .
	\end{split}
\end{equation}
If we combine this with the straightforward estimates
\begin{align*}
  \| \skew \sigma \|_{L^2 (\Omega)} & = \| \skew (\dev\sigma - \mu \sym{\rm grad\,} u) \|_{L^2 (\Omega)}
  \leq \| \dev\sigma - \mu \sym{\rm grad\,} u \|_{L^2 (\Omega)} \: , \\
  \mu \| {\rm div} \: u \|_{L^2 (\Omega)} & = \| {\rm tr} (\dev\sigma - \mu\sym {\rm grad\,} u) \|_{L^2 (\Omega)}
  \leq \sqrt{3} \| \dev\sigma - \mu\sym {\rm grad\,} u \|_{L^2 (\Omega)} \: ,
\end{align*}
we are led to
\begin{align*}
	\| & \dev\sigma-\mu\sym{\rm grad\,} u \|_{L^2(\Omega)}^2 + \|\Div\sigma\|_{L^2(\Omega)}^2 \\
	& \geq \frac{1}{3} \left( \| \dev\sigma-\mu\sym{\rm grad\,} u \|_{L^2(\Omega)}^2
	+ \| \skew \: \sigma \|_{L^2(\Omega)}^2 + \frac{\mu^2}{3} \| {\rm div} \: u \|_{L^2(\Omega)}^2
	+ \|\Div\sigma\|_{L^2(\Omega)}^2 \right) \\
	& \geq \frac{\delta}{6} \left( \| \dev\sigma-\mu\sym{\rm grad\,} u \|_{L^2(\Omega)}^2
	+ \frac{2}{\delta} \left( \| \skew \: \sigma \|_{L^2(\Omega)}^2
	+ \frac{\mu^2}{3} \| {\rm div} \: u \|_{L^2(\Omega)}^2
	+ \| \Div\sigma \|_{L^2(\Omega)}^2 \right) \right) \\
         & \geq \frac{\delta}{6} \left( \| \dev\sigma \|_{L^2(\Omega)}^2
         + \mu^2 \| \sym{\rm grad\,} u \|_{L^2(\Omega)}^2
         - \delta \left( \mu^2 \| u \|_{L^2 (\Omega)}^2 + \mu^2 \| {\rm grad\,} u \|_{L^2 (\Omega)}^2
                                                    + \frac{1}{3} \| \tr \sigma \|_{L^2 (\Omega)}^2 \right) \right. \\
         & \qquad \left. + \frac{1}{\delta} \left( \| \Div \sigma \|_{L^2 (\Omega)}^2
                                                    + \| \skew \sigma \|_{L^2 (\Omega)}^2
                                                    + \frac{1}{3} \mu^2 \| {\rm div} \: u \|_{L^2 (\Omega)}^2 \right) \right) \: \\
         & \geq \frac{\delta}{6} 
         \left( 
         \| \dev\sigma \|_{L^2(\Omega)}^2
         + \| \Div\sigma \|_{L^2(\Omega)}^2
         + \mu^2 \| \sym{\rm grad\,} u \|_{L^2(\Omega)}^2
         - \delta \mu^2 \| u \|_{H(\Grad;\Omega)}^2 
         - \frac{\delta}{3} \| \tr \sigma \|_{L^2 (\Omega)}^2 
         \right) \\
         & \geq \frac{\delta}{6} 
         \left( 
         \frac{1}{C_{D\!D}^2}\| \sigma \|_{H(\Div;\Omega)}^2
         + \frac{\mu^2}{C_{K\!P}^2} \| u \|_{H(\Grad;\Omega)}^2
         - \delta \mu^2 \| u \|_{H(\Grad;\Omega)}^2 
         - \delta \|  \sigma \|_{L^2 (\Omega)}^2 
         \right)
\end{align*}
for all $\delta\leq1$ with $C_{D\!D}$ from Theorem \ref{Dev-Div-T} and the Korn-Poincar\'e constant $C_{K\!P}$
in the Korn-Poincar\'e inequality 
$$C_{K\!P}\| \sym{\rm grad\,} u \|_{L^2(\Omega)}\geq C_{K} \| {\rm grad\,} u \|_{L^2(\Omega)}\geq  \|  u \|_{H(\Grad;\Omega)}.$$		
Choosing $\delta$ sufficiently small gives us the desired coercivity estimate
\begin{equation}
  \| \dev\sigma-\mu\sym{\rm grad\,} u \|_{L^2(\Omega)}^2 + \| \Div\sigma \|_{L^2(\Omega)}^2
   \geq C
  \left( \| \Div\sigma \|_{H(\Div;\Omega)}^2 
  + \| u \|_{H (\Grad;\Omega)}^2 \right) \: .
  \label{eq:coercivity_pseudostress}
\end{equation}
The pseudostress-velocity formulation of the stationary Stokes equations introduced above was
studied in \cite[section 3.2]{caileewang:04} (see also \cite{caitongvassilevskiwang:10} and
\cite{gaticamarquezsanchez:10} for related mixed finite element approaches). It was used as a basis
for the treatment of Stokes-Darcy interface problems by a first-order system least squares approach
in \cite{muenzenmaierstarke:11}. Recently, a pseudostress-based approach for the stationary
Navier-Stokes was investigated in \cite{caizhang:12}.

\subsection{Pseudostress formulation of generalized Newtonian flow}

The estimate of Theorem \ref{Dev-Div-T} is also useful in the context of nonlinear generalized
Newtonian fluids which differs from the formulation above in that the viscosity may depend on the velocity
gradient $\mu = \mu ({\rm grad\,} u)$. Very popular is Carreau's law, where this nonlinear dependence is
given by
\[
  \mu ({\rm grad\,} u) = \mu_0 \left( 1 + | \sym {\rm grad\,} u |^2 \right)^{(r-2)/2}
\]
with $\mu_0 > 0$ and $r \geq 1$. Depending on the value of $r$, shear-thickening or shear-thinning
behavior of the fluid can be modeled. A dual-mixed approach to nonlinear generalized Newtonian
Stokes flow was introduced and analyzed in \cite{ervinhowellstanculescu:08}. This model may also
be treated by a pseudostress-velocity approach in a first-order system least squares setting based on
minimizing the nonlinear functional
\begin{equation}
  \|\dev\sigma-\mu({\rm grad\,} u) \sym{\rm grad\,} u\|_{L^2(\Omega)}^2 + \|\Div\sigma-f\|_{L^2(\Omega)}^2 \: .
\end{equation}
Such a method is studied in detail in \cite{muenzenmaier:12}.

\subsection{Infinitesimal gradient plasticity}

Phenomenological plasticity models are intended to describe the irreversible deformation behavior of metals. 
There exists a great variety of models. Here we focus on rate-dependent or rate-independent models with kinematic hardening. The system of equations consist of balance of linear momentum coupled with a local nonlinear evolution equation in each space point for the plastic variable.

In many new applications, the size of the considered specimen is so small, that size effects need to be taken into account. Instead of a local evolution problem we have to consider a nonlinear evolution problem where the right hand side contains certain combinations of second partial derivatives of the plastic distortion. 

For the setting of the nonlinear gradient-plasticity problem,
let $\Omega \subset \R^3$ be an open and bounded set, the set of material
points of the solid body. By $T_e$ we denote a positive number (time of existence). 
Unknown in our small strain formulation are the displacement field $u:\Omega\times [0,T_e)\to \R^3$ of the material point $x$ at time $t$ and the non-symmetric infinitesimal
plastic distortion $P: \Omega\times [0,T_e)\to  \mathfrak{sl}(3)$. The model equations of the problem are
\begin{eqnarray}
 \Div \sigma &=&  f, \notag\\
\sigma &=& 2\mu\, ( \sym ({\rm grad\,} u - P ) )+\lambda\, \tr{{\rm grad\,} u}\cdot{\rm Id} ,   \notag\\
\label{CurlPr3} 
\partial_t P(x,t) & \in & g \big(x,\Sigma^{\rm lin}(x,t)\big),  \\
 \Sigma^{\rm lin}& = &\Sigma^{\rm lin}_{\rm e}+\Sigma^{\rm lin}_{\rm sh}
+\Sigma^{\rm lin}_{\rm curl},\label{microPr3} \notag\\
\Sigma^{\rm lin}_{\rm e} &= &\sigma, \quad
\Sigma^{\rm lin}_{\rm sh}=- \dev \sym P,\quad
\Sigma^{\rm lin}_{\rm curl}=- \Curl\Curl P\, , \notag
\end{eqnarray}
which must be satisfied in $\Omega \times [0,T_e)$. Here, $\Sigma^{\rm lin}$ is the infinitesimal Eshelby stress tensor driving the evolution of the plastic distortion $P$. The initial and boundary conditions are  
\begin{align*} 
P(\,\cdot\,,0) &= P_{0} && \text{in }\Omega, \\
\nu\times P &=0&&\text{on }\partial\Omega \times [0,T_e), \label{CurlPr5}\\
u &= 0 && \text{on }\partial \Omega \times  [0,T_e)\,,\notag
\end{align*}
where $\nu$ is a normal vector on the boundary $\partial\Omega$. For the model we require that the nonlinear constitutive mapping
$(\Sigma\to g(\,\cdot\,, \Sigma)):\R^{3\times 3} \rightarrow 2^{\mathfrak{sl}(3)}$ is monotone. Given are the volume force 
$f:\Omega\times [0,T_e)\to \R^3$ and the initial datum $P_0 :\Omega\to\mathfrak{sl}(3)$.
It is easy to see that the corresponding free energy of the system is
\begin{align}
  \mu \|\sym ({\rm grad\,} u -P)\|_{L^2(\Omega)}^2+\frac{\lambda}{2}\|\tr{{\rm grad\,} u}\|_{L^2(\Omega)}^2
  +\frac{1}{2} \|\dev\sym P\|_{L^2(\Omega)}^2+\frac{1}{2} \|\Curl P\|_{L^2(\Omega)}^2\, .
\end{align}
The appearance of $\Curl P$ instead of the full gradient ${\rm grad\,} P$ is dictated by dislocation mechanics, the appearance of 
$\dev\sym P $ instead of $P$ is dictated by invariance of the model under superposition of infinitesimal rotations. Here, coercivcity is obtained by using the DevSym-Curl inequality. Model equations similar to the above problem have 
been considered in \cite{Nesenenko_Neff12a,Nesenenko_Neff12b,Ebobisse_Neff09,Neff_Sydow_Wieners08,Neff_Chelminski07_disloc}.

\subsection{Infinitesimal Cosserat elasticity}

Cosserat or micropolar elasticity is intended to describe materials with a microstructure which has the degrees of freedom of a rigid body. With Cosserat elasticity, it is possible to describe some form of elastic size effects (smaller samples are comparatively stiffer) and wave dispersion in the case of dynamic equations. Here, we consider the static problem. In a variational context, the problem is completely described by writing the energy which is to be minimized. We are looking for the displacement $u:\Omega\to\R^3$ and the infinitesimal Cosserat microrotation $A:\Omega\to\so(3)$ minimizing the two-field functional
\begin{align*}
  \mu \|\sym {\rm grad\,} u\|_{L^2(\Omega)}^2
  +\mu_{c} \|\skew ({\rm grad\,} u -A)\|_{L^2(\Omega)}^2
  +\frac{\lambda}{2}\|\tr{{\rm grad\,} u}\|_{L^2(\Omega)}^2
  +\frac{1}{2} \|\dev\sym\Curl A\|_{L^2(\Omega)}^2
  +(f,u)_{L^2(\Omega)}\, .
\end{align*}
The corresponding system of Euler-Lagrange equations in strong form are
\begin{align*}
   \Div \sigma&=f\,, \\
\sigma &= 2\mu\, \sym {\rm grad\,} u+\lambda\, \tr{{\rm grad\,} u}\cdot{\rm Id} +2\mu_c\skew({\rm grad\,} u-A),   \\
   \skew\sigma&=\skew\Curl\dev\sym \Curl A\, .
\end{align*}
The form of the curvature contribution $\dev\sym\Curl A$ instead of the full gradient ${\rm grad\,} A$ is motivated by conformal invariance of the model, see \cite{Neff_Jeong_Conformal_ZAMM08}. Here, a variant of the DevSym-DevCurl inequality is applicable. Model equations similar to the above problem have 
been considered in \cite{Neff_Jeong_Conformal_ZAMM08,Neff_JeongMMS08,Neff_ZAMM05}.

\subsection{Infinitesimal Cosserat elasto-plasticity}

Frequently encountered  are also couplings between Cosserat elasticity and plasticity models. However, plasticity in these models is treated classically as a local phenomenon. 
We are looking for the displacement $u:\Omega\times [0,T_e)\to\R^3$, the infinitesimal Cosserat micro-rotation
$A:\Omega\times [0,T_e)\to\so(3)$ and the plastic distortion $P:\Omega\times [0,T_e)\to\sL(3)$ satisfying 
\begin{eqnarray}
\Div \sigma &=&  f, \notag  \\
\sigma &=& 2\mu\, ( \sym ({\rm grad\,} u - P ) )+\lambda\, \tr{{\rm grad\,} u}\cdot{\rm Id}
 +2\mu_c\skew({\rm grad\,} u-A),   \notag  \\
   \partial_t P(x,t) & \in & g \big(x,\sym\sigma(x,t)\big),   \\
  \skew\sigma&=&\skew\Curl\dev \Curl A\, .   \notag
  \end{eqnarray}
Model equations with these features have 
been considered in \cite{Neff_Knees06,Neff_Chelminski03a,Neff_Cosserat_plasticity05,Neff_Chelminski_qam07,Neff_Chelminski_Wieners06} with the purpose of obtaining regularizations of classical plasticity models.

\subsection{Infinitesimal relaxed micromorphic elasticity}

Micromorphic extended continuum models assume that at each material point there is a microstructure attached which itself may deform as an elastic body. 
In a variational context, the problem is completely described by writing down the energy which is to be minimized. We are looking for the displacement fields $u:\Omega\to\R^3$ and the not necessarily symmetric micromorphic distortion $P:\Omega\to\R^{3\times 3}$ minimizing
\begin{align*}
  \mu \|\sym ({\rm grad\,} u-P)\|_{L^2(\Omega)}^2
 +\frac{\lambda}{2}\|\tr{{\rm grad\,} u-P}\|_{L^2(\Omega)}^2
 +\frac{1}{2} \|\dev\sym P\|_{L^2(\Omega)}^2
  +\frac{1}{2} \|\Curl P\|_{L^2(\Omega)}^2
  +(f,u)_{L^2(\Omega)}\, .
  \end{align*}
The corresponding system of Euler-Lagrange equations in strong form are
\begin{align*}
   \Div \sigma&=f\,, \\
\sigma& = 2\mu\, \sym ({\rm grad\,} u-P)+\lambda\, \tr{{\rm grad\,} u-P}\cdot{\rm Id},    \\
   \Curl\dev \Curl P&=-\dev\sym P+ \sigma\, .
   \end{align*}
An important feature, which sets this model apart from more classical micromorphic approaches, is that the balance of forces does not `see' derivatives of $P$ since $\Div\Curl=0$. Here, the DevSym-DevCurl-inequality is applicable. Model equations similar to the above problem have 
been considered \cite{Neff_Jeong_ZAMP08,Neff_micromorphic_rse_05,GhibaNeff:13}.

\appendix

\section{The kernel of dev\,sym\,Grad}

	For the convenience of the reader we compute the representation formulae (\ref{RepresentationFormulaFunction_T})
	and (\ref{RepresentationpA_T}) of vector fields in the kernel of ${\rm dev\,}{\rm sym\;}{\rm Grad}$,
	 used in the proof of Lemma \ref{DevSymGrad-T}. These mappings are often called conformal mappings or 
	conformal Killing vectors. For $n=3$ such a representation is given, e.g. in \cite{Neff_JeongMMS08} and for 
	arbitrary $n\geq 3$ in \cite{schirra:12}.
	 Now let $\Omega\subset \R^n$ be an arbitrary domain, $n\geq 3$ and 
 	$v\in H({\rm Grad\,};\Omega)$ with ${\rm dev\,}{\rm sym\;}{\rm Grad\,}v=0$.	
	 Then
		\begin{equation}\label{RepresentationNablaPhi}
			{\rm Grad\,}v = \p\cdot{\rm Id}+A\,, 
		\end{equation}
	where after selecting a suitable representant
	$\p(x)\in\R$  and $A(x)$ is  a skew-symmetric $(n\times n)$-matrix for all $x\in \Omega$.
	 Since ${\rm Curl\,}{\rm Grad\,}v=0$ we obtain that for all $i,j,k=1, \dots n$
		\begin{equation}\label{CurlZero}
			\left(\partial_j \p\right)\delta_{ik}-\left(\partial_k \p\right)\delta_{ij} = -\partial_j A_{ik}+\partial_k A_{ij}\,.
		\end{equation}
	Now assume that $i,j,k$ are pairwise different, then using (\ref{CurlZero}) and $A_{ij}=-A_{ji}$ we compute
		\begin{equation}\nonumber
			\partial_j A_{ik}=\partial_k A_{ij}=-\partial_k A_{ji}=-\partial_i A_{jk}=\partial_i A_{kj}=\partial_j A_{ki}=-\partial_j A_{ik} \,,
			\end{equation}
	yielding $\partial_j A_{ik}=0$. Now assume $j=i$, but $k\not=i$. In this case we obtain by \eqref{CurlZero}
		\begin{equation}\nonumber
			\partial_j A_{jk}= \partial_k \p\,.
		\end{equation}
Therefore,
		\begin{equation}\label{RelationsAij}
			\begin{array}{rclcl}
				\partial_j A_{jk} & = &	\partial_k \p &\text{if}  &  j\not= k,\\ 
				\partial_j A_{kj} & = & -\partial_k \p &\text{if} &  j\not= k,\\ 
				\partial_j A_{ik} & = & 0 & \text{if}& i\not = j, i\not= k, j\not=k,\\ 
				\partial_j A_{ii}	 & = & 0\,. 
			\end{array}
		\end{equation}
	Now we show that (\ref{RelationsAij}) implies ${\rm grad\,} \p = const$.
	First assume $j\not=k$ and choose $i$ with $i\not=j$ and $i\not=k$ (since $n\geq 3$, this is possible). 
	Therefore we obtain
		\begin{equation}\nonumber
			\partial_k\partial_j \p = \partial_k\partial_i A_{ij}= \partial_i\partial_k A_{ij}=0\,.
		\end{equation}
	Now, we assume that $i\not=j$, then
		\begin{equation}\nonumber
			\partial_j\partial_j \p = \partial_j\partial_i A_{ij}=-\partial_j\partial_i A_{ji}
		=-\partial_i\partial_j A_{ji}=-\partial_i\partial_i \p\,.
		\end{equation}
	As $n\geq 3$ we can play the indices against each other and obtain 
		$$\partial_j\partial_j \p=-\partial_i\partial_i \p=\partial_k\partial_k \p=-\partial_j\partial_j \p=0\,$$
		for $i,j,k$ pairwise different.
	Therefore ${\rm grad\,} \p = const = \bar{w}\in \R^n$ and after a possible redefinition
	on a set of measure zero, we get
		\begin{equation}\label{Representationp}
			\p(x) = \bar{w}\cdot x +\bar{\p}
		\end{equation}
	with $\bar{\p}\in\R$. Note $\partial_{i}\p=\bar{w}_{i}$.
	Since $\partial_k A_{ij}$ is constant, see (\ref{RelationsAij}), we also know that 
		\begin{equation}\label{RepresentationA}
			A_{ij}(x) = \sum_{k=1}^n\partial_k A_{ij}x_k+\bar{A}_{ij}
			=\sum_{k=1}^n \bar{a}_{ijk}x_k+\bar{A}_{ij}
			=\bar{w}_{j}x_{i}-\bar{w}_{i}x_{j}+\bar{A}_{ij}
		\end{equation}		
	with $\bar{a}_{ijk}$ from \eqref{aijk} and some skew-symmetric constant  matrix $\bar{A}$.	
	Utilizing (\ref{RepresentationNablaPhi}), (\ref{Representationp}) and
	 (\ref{RepresentationA}) we conclude that (\ref{RepresentationFormulaGradient_T}) holds true. 
	 Furthermore, by integrating the $i$-th component of $v$ 
	 we obtain from \eqref{RepresentationNablaPhi}, i.e.,
	$$\partial_{j}v_{i}(x)=u(x)\delta_{ij}+A_{ij}(x)=u(x)\delta_{ij}+\bar{w}_{j}x_{i}-\bar{w}_{i}x_{j}+\bar{A}_{ij},$$
	immediately
  $$v_{i}(x) = u(x)x_{i}-\frac{1}{2}\bar{w}_{i}|x|^2+\bar{A}_{ik}x_{k}+\bar{v}_{i}$$
  or as a vector
  $$v(x) = u(x)x-\frac{1}{2}|x|^2\bar{w}+\bar{A}x+\bar{v}$$
which is \eqref{RepresentationFormulaFunction_T}.  
%

{\footnotesize
\bibliographystyle{plain} 
\bibliography{patrizio,sebastian}}

\end{document}